\documentclass[preprint,12pt]{elsarticle}

\usepackage{amssymb}
\usepackage{amsmath}
\usepackage{amsthm}
\usepackage{enumitem}
\let\ker\undefined
\DeclareMathOperator{\ker}{ker}

\renewcommand\div{\operatorname{div}}

\newcommand\curl{\operatorname{curl}}

\newcommand\rot{\operatorname{rot}}

\newcommand\Rey{\operatorname{Re}}
\newcommand\Rem{\operatorname{Re}_m}

\newcommand{\dx}{\,\mathrm{d}x}

\usepackage{bm}
\usepackage{tikz-cd}
\usepackage{tikz}
\usepackage{forest}
\usepackage{pgfplots}
\pgfplotsset{compat=1.18}
\usetikzlibrary{trees,positioning,arrows.meta}
\usepgfplotslibrary{groupplots}

\tikzset{box/.style={draw,align=left,inner sep=4pt,fill=white},pink/.style={box,fill=red!12},blue/.style={box,fill=blue!14}}
\usepackage{subcaption}
\usepackage[ruled]{algorithm}

\newtheorem{theorem}{Theorem}[section]
\newtheorem{proposition}{Proposition}[section]

\newtheorem{remark}{Remark}[section]
\newtheorem*{remark*}{Remark}
\newtheorem{lemma}[theorem]{Lemma}
\newcommand{\ov}[1]{\overline{#1}}
\usepackage{mathtools}
\usepackage{empheq}
\usepackage{booktabs}
\usepackage{multirow}
\usepackage{siunitx}
\usepackage{hyperref}
\hypersetup{hidelinks}
\usepackage[nameinlink]{cleveref}
\crefname{algorithm}{Algorithm}{Algorithms}
\Crefname{algorithm}{Algorithm}{Algorithms}
\usepackage{todonotes}

\usepackage[breakable,skins]{tcolorbox}
\newif\ifshowannotations
\showannotationstrue
\newtcolorbox{MDHbox}{breakable,enhanced,colback=orange!8,colframe=orange!75!black,
  title=MDH,fonttitle=\bfseries\small,arc=2pt,boxrule=0.8pt,
  left=4pt,right=4pt,top=3pt,bottom=3pt,before skip=6pt,after skip=6pt}

\journal{Journal of Computational Physics}

\begin{document}

\begin{frontmatter}



\title{Structure-preserving finite element method for the 2D MHD-$\alpha$ turbulence system}

\author[label1,label2]{Patrick E.~Farrell}
\author[label1]{Mingdong He}
\author[label1]{Kaibo Hu}
\affiliation[label1]{organization={Mathematical Institute, University of Oxford},
            addressline={Andrew Wiles Building, Radcliffe Observatory Quarter},
            city={Oxford},
            postcode={OX2 6GG},
            country={UK}}

\affiliation[label2]{organization={Mathematical Institute, Faculty of Mathematics and Physics, Charles University},
            city={Prague},
            country={Czechia}}

\begin{abstract}
We propose a structure-preserving finite element method for the two-dimensional MHD-$\alpha$ turbulence model, which filters the transport velocity at a length $\alpha$ while retaining the ideal invariants of the original MHD equations, with the energy replaced by an $\alpha$-dependent energy. The scheme preserves mass conservation and
the magnetic Gauss law pointwise and satisfies discrete dissipation laws for the energy, the
cross helicity and the mean square magnetic potential, which reduce to exact conservation in
the ideal limit. The discrete systems are solved by an augmented Lagrangian preconditioner
with monolithic geometric multigrid, observed to be robust with respect
to the filter width $\alpha$ and the Reynolds numbers, and to the mesh size for small time steps. The scheme converges at the optimal rates and agrees
with a Fourier pseudo-spectral direct simulation of the original MHD equations at the large scales. In under-resolved flows, the unfiltered scheme has no numerical dissipation and accumulates energy at the grid scale. Filtering relieves this accumulation and a larger $\alpha$ gives a more compact spectrum. Long-time simulations reproduce selective
decay, including relaxation toward the Arnold-type bound and the inverse cascade of the mean
square potential. 

\end{abstract}




\end{frontmatter}
\section{Introduction}
The long-time dynamics of magnetohydrodynamic (MHD) turbulence is governed by its ideal
invariants \cite{biskamp2003magnetohydrodynamic}. In three dimensions these are the total
energy, the cross helicity, and the magnetic helicity. In two dimensions the magnetic helicity
vanishes, and its role is taken by the mean square magnetic potential. We consider incompressible
MHD in two dimensions on the periodic domain $\Omega=\mathbb T^2$, written in terms of the velocity $\bm u$ and the magnetic potential $A$,
\begin{subequations}\label{eqn:mhd}
\begin{align}
&\partial_t\bm u+(\bm u\cdot\nabla)\bm u+\nabla p=S(\bm B\cdot\nabla)\bm B+\nu\Delta\bm u+\bm f,\\
&\partial_t A+\bm u\cdot\nabla A=\eta\Delta A,
\end{align}
\end{subequations}
with $\nabla\cdot\bm u=0$ and magnetic field $\bm B=\curl A$. Here $\nu=\Rey^{-1}\ge0$ and $\eta=\Rem^{-1}\ge0$ are the
nondimensional viscosity and resistivity, the inverse kinetic and magnetic Reynolds numbers, $S>0$ is the coupling coefficient and $\bm f$ is a prescribed forcing term.
We use the two-dimensional curl operators $\curl\phi=(\partial_y\phi,-\partial_x\phi)^T$ and
$\rot\bm v=\partial_xv_2-\partial_yv_1$, and denote by $w=\rot\bm u$ the vorticity and by
$j=\rot\bm B=-\Delta A$ the current density. For $\bm a,\bm b\in\mathbb R^2$ and $\phi\in\mathbb R$ we write
$\bm a\times\phi=-\phi\times\bm a=(a_2\phi,-a_1\phi)^T$ and $\bm a\times\bm b=a_1b_2-a_2b_1$.
On $\mathbb T^2$ a divergence-free field is the curl of a periodic potential if and only if its
mean vanishes, the constant vector fields being the harmonic fields of the torus. We therefore
assume that $\bm u$, $\bm B$ and $\bm f$ have zero mean. Since $\bm B=\curl A$ determines $A$ only up to a constant, we fix it by $\int_\Omega A\dx=0$. For
$\nu=\eta=0$ and $\bm f=0$, \eqref{eqn:mhd} conserves the energy, the cross helicity and the mean
square potential,
\begin{equation}\label{eqn:mhd-invariants}
\mathcal E=\frac12\int_\Omega|\bm u|^2+S|\bm B|^2\dx,\qquad
\mathcal H_c=\int_\Omega\bm u\cdot\bm B\dx,\qquad
\mathcal A=\int_\Omega A^2\dx.
\end{equation}
With small viscosity and resistivity, these three invariants have been observed to decay at different rates in decaying turbulence \cite{Matthaeus1980,biskamp1994dynamics,biskamp2003magnetohydrodynamic}; the system then relaxes
toward a state of minimal energy constrained by them, as observed in simulations of
two-dimensional MHD \cite{biskamp1989dynamics}. This so-called selective decay mechanism underlies the inverse cascade
of magnetic helicity in three-dimensional MHD turbulence \cite{frisch1975possibility} and of the
mean square potential in two-dimensional MHD turbulence
\cite{pouquet1978two,fyfe1976high,biskamp1989dynamics}, and the relaxation toward Taylor states
\cite{taylor1986relaxation} in fusion plasmas.

In a turbulent flow, energy cascades from the large scales toward ever smaller scales until it
is removed by viscosity at the dissipation scale according to Kolmogorov's 1941 theory \cite{kolmogorov1991local}. A direct numerical
simulation (DNS) resolves all scales down to this dissipation scale, so its cost grows rapidly
with the Reynolds number and becomes prohibitive \cite{pope2000turbulent}. Large eddy simulation (LES) instead resolves only the large scales and
models the effect of the unresolved small scales on them \cite{sagaut2006large}. The
fundamental requirement is to truncate the scales to those representable on
a computationally feasible mesh without substantially altering the large
flow structures~\cite{layton2012approximate}.

For long-time turbulence simulations to be reliable, the numerical method should respect the
invariants and their dissipation laws, which calls for structure-preserving discretizations.
A structure-preserving scheme, however, introduces no numerical dissipation. While the flow is
resolved, this is an asset; but in long-time turbulent computations the energy it conserves cannot be removed there and accumulates at the smallest resolved scales, thus polluting the simulation. Conventional turbulence codes address under-resolution at the level of the
scheme, through numerical or eddy viscosity, which, however, changes the dissipation rates of the ideal invariants. This highlights an interesting tension between structure-preserving schemes and the need for dissipation in under-resolved turbulence simulations. To address this issue, we seek to suppress the accumulation of energy at the smallest resolved scales without introducing additional dissipation, while retaining the cross helicity and the mean square potential~\eqref{eqn:mhd-invariants} of the original MHD system, together with a modified energy, so that long-time computations remain reliable. This can be achieved with the regularized $\alpha$-model \cite{linshiz2007analytical,foias2001navier}, which we now review from the perspective of the ideal MHD invariants.

We aim to regularize the system \eqref{eqn:mhd} via filtering the transport velocity. The simplest way is to apply a Helmholtz filter to the velocity
\begin{equation}\label{eqn:filter}
\ov{\bm u}:=(1-\alpha^2\Delta)^{-1}\bm u,
\end{equation}
which is divergence free. In the ideal case, the potential is transported by $\ov{\bm u}$,
\begin{equation}\label{eqn:A-transport}
\partial_tA+\ov{\bm u}\cdot\nabla A=0,
\end{equation}
and we write the momentum equation as $\partial_t\bm u+\bm N+\nabla p=S(\bm B\cdot\nabla)\bm B$
with an unknown transport term $\bm N$. Our goal is to derive a transport term $\bm N$ such that the three invariants \eqref{eqn:mhd-invariants} are preserved like the original MHD. 

\emph{Mean square potential.} Testing \eqref{eqn:A-transport} with $2A$ and using
$\nabla\cdot\ov{\bm u}=0$,
\begin{equation}
\frac{d}{dt}\mathcal A=-\int_\Omega\ov{\bm u}\cdot\nabla(A^2)\dx=\int_\Omega(\nabla\cdot\ov{\bm u})\,A^2\dx=0,
\end{equation}
for every $\bm N$. Transporting $A$ by any divergence-free velocity conserves $\mathcal A$, and its
dissipation law is unchanged when $\eta\Delta A$ is added. 

\emph{Cross helicity.} It is easy to derive that
\begin{equation}
\frac{d}{dt}\mathcal H_c=\int_\Omega\partial_t\bm u\cdot\bm B+w\,\partial_tA\dx=-\int_\Omega\bm B\cdot\big(\bm N+\ov{\bm u}\times w\big)\dx.
\end{equation}
For this to vanish for every divergence-free $\bm B$, the bracket must be a gradient, absorbed into the pressure, plus a constant vector, which vanishes since $\bm u$ keeps zero mean. Hence
\begin{equation}\label{eq:N}
\bm N=(\ov{\bm u}\cdot\nabla)\bm u+\sum_{j=1}^2\bm u_j\nabla\ov{\bm u}_j
=-\ov{\bm u}\times w+\nabla(\ov{\bm u}\cdot\bm u),
\end{equation}
which reduces to $(\bm u\cdot\nabla)\bm u$ up to a gradient when $\alpha=0$. Thus, we obtain the transport term. 

\emph{Energy.} By \eqref{eqn:A-transport} and $(\bm B\cdot\nabla)\bm B=j\times\bm B+\nabla(|\bm B|^2/2)$ with $\ov{\bm u}\cdot(j\times\bm B)=j\,\ov{\bm u}\cdot\nabla A$, the magnetic energy changes by the work of the Lorentz force on $\ov{\bm u}$,
\begin{equation}
\frac S2\frac{d}{dt}\int_\Omega|\bm B|^2\dx=S\int_\Omega j\,\partial_tA\dx=-S\int_\Omega\ov{\bm u}\cdot(\bm B\cdot\nabla)\bm B\dx.
\end{equation}
We therefore test the momentum equation with $\ov{\bm u}$: the transport term \eqref{eq:N} does no work since $(\ov{\bm u}\times w)\cdot\ov{\bm u}=0$, the Lorentz terms cancel, and by the symmetry of the filter the conserved energy is
\begin{equation}\label{eqn:E-alpha}
\mathcal E^\alpha=\frac12\int_\Omega\ov{\bm u}\cdot\bm u+S|\bm B|^2\dx
=\frac12\int_\Omega|\ov{\bm u}|^2+\alpha^2|\nabla\ov{\bm u}|^2+S|\bm B|^2\dx.
\end{equation}
Adding viscosity, resistivity and forcing, we obtain the MHD-$\alpha$ model
\cite{linshiz2007analytical}
\begin{subequations}\label{eqn:mhd-alpha}
\begin{align}
&\partial_t\bm u+(\ov{\bm u}\cdot\nabla)\bm u+\sum_{j=1}^2\bm u_j\nabla\ov{\bm u}_j+\nabla p
=S(\bm B\cdot\nabla)\bm B+\nu\Delta\bm u+\bm f,\label{eqn:mhd-alpha-u}\\
&\partial_tA+\ov{\bm u}\cdot\nabla A=\eta\Delta A,\label{eqn:mhd-alpha-A}
\end{align}
\end{subequations}
with $\ov{\bm u}=(1-\alpha^2\Delta)^{-1}\bm u$, $\nabla\cdot\bm u=0$ and $\bm B=\curl A$. It is the
magnetohydrodynamic analog of the Navier--Stokes-$\alpha$ model
\cite{chen1998camassa,foias2001navier}.

\begin{remark*}[Kelvin circulation]
The transport term \eqref{eq:N} is the one produced by Kelvin's circulation theorem when material
loops are advected by $\ov{\bm u}$ instead of $\bm u$ \cite{foias2001navier,linshiz2007analytical}.
A line element joining two such particles evolves by
$\frac{d}{dt}(d\bm x)=\ov{\bm u}(\bm x+d\bm x)-\ov{\bm u}(\bm x)=(d\bm x\cdot\nabla)\ov{\bm u}$, so
for a loop $\gamma(t)$ moving with $\ov{\bm u}$,
\begin{equation}
\begin{aligned}
\frac{d}{dt}\oint_{\gamma(t)}\bm u\cdot d\bm x
&=\oint_{\gamma(t)}\big(\partial_t\bm u+(\ov{\bm u}\cdot\nabla)\bm u\big)\cdot d\bm x+\oint_{\gamma(t)}\bm u\cdot\frac{d}{dt}(d\bm x)\\
&=\oint_{\gamma(t)}\big(\partial_t\bm u+(\ov{\bm u}\cdot\nabla)\bm u\big)\cdot d\bm x+\oint_{\gamma(t)}\bm u\cdot(d\bm x\cdot\nabla)\ov{\bm u}\\
&=\oint_{\gamma(t)}\Big(\partial_t\bm u+(\ov{\bm u}\cdot\nabla)\bm u+\sum_{j=1}^2\bm u_j\nabla\ov{\bm u}_j\Big)\cdot d\bm x.
\end{aligned}
\end{equation}
For $\alpha=0$ the last term is $\nabla(|\bm u|^2/2)$, which integrates to zero around the loop, and
the classical theorem is recovered.
\end{remark*}

The MHD-$\alpha$ model \eqref{eqn:mhd-alpha} thus retains the three ideal invariants of MHD, with $\mathcal E$ replaced by
$\mathcal E^\alpha$. The term $\alpha^2|\nabla\ov{\bm u}|^2$ in $\mathcal E^\alpha$ penalizes
motions below the length $\alpha$, so small-scale velocity activity below this scale is suppressed,
which may reduce the resolution required \cite{foias2001navier,chen1999direct}, while no dissipation
is added. The cross helicity $\mathcal H_c$ and the mean square potential $\mathcal A$ are unchanged
and satisfy the same dissipation laws as in MHD.

To capture these dynamics, a discretization should preserve the ideal invariants and their
nonideal dissipation laws at the discrete level, since numerical dissipation acting on the slowly decaying invariants alters the very hierarchy that drives selective decay. The past decades have seen significant progress in structure-preserving methods for flows and MHD. In particular, schemes based on the finite element exterior calculus (FEEC)
\cite{arnoldFiniteElementExterior2006,arnoldFiniteElementExterior2010,ArnoldFiniteElementExterior2018}
have been developed that precisely preserve important structure, such as the magnetic Gauss law
and helicity conservation \cite{huStableFiniteElement2017,hiptmair2018fully,huHelicityconservativeFiniteElement2021,gawlikFiniteElementMethod2022,LaakmannStructurepreservinghelicityconservingfinite2023,maoIncompressibilityDivB0Preserving2025,da2025error,li2026linear}.
Helicity-preserving discretizations have been used to study magnetic relaxation
\cite{he2025helicity,parkerE3}, where preserving the magnetic helicity at the discrete level
is essential: the ideal invariant constrains the evolution and determines the relaxed state.
The same method has also been used to compute three-dimensional MHD equilibria of
stellarators \cite{BlickhanMRXdifferentiable3D2025}. Energy- and helicity-conserving schemes have also been constructed with mixed finite elements for the
Navier--Stokes equations \cite{RebholzEnergyHelicityConservingFinite2007} and the Navier--Stokes-$\alpha$ model
\cite{miles2010enhanced}, and for the Navier--Stokes equations with mimetic dual-field discretizations
\cite{zhangMassKineticEnergy2022,zhangMEEVCDiscretizationTwodimensional2024}. Beyond finite elements, long-time
simulations of ideal MHD turbulence were carried out in \cite{modin2026structure-mhd} with a
structure-preserving matrix discretization, and a structure-preserving spectral scheme for ideal MHD
was proposed in \cite{guo2026novel}. 

In this work we construct a structure-preserving finite element method for the two-dimensional
MHD-$\alpha$ model \eqref{eqn:mhd-alpha}, and show that both the filter and structure preservation are needed for accurate simulations on under-resolved meshes, i.e.\ for use as a large eddy simulation. Our contributions are:
\begin{itemize}
\item \emph{The discretization.} A structure-preserving finite element discretization of the
  MHD-$\alpha$ model with pointwise divergence-free fields and discrete
  dissipation laws for all three ideal invariants, which reduce to exact conservation in the ideal
  limit.
\item \emph{The solver.} An augmented Lagrangian preconditioner with a monolithic multigrid solver, whose iteration counts are observed to be robust with respect to the filter width and the Reynolds numbers, and to the mesh size for small time steps.
\item \emph{Why the MHD-$\alpha$ model is necessary.} A comparison with a direct numerical
  simulation shows that, once the flow is under-resolved, the unfiltered scheme accumulates energy
  at the grid scale, which pollutes the small scales with spurious oscillations. The filter suppresses
  this accumulation without introducing additional dissipation and with all three invariants preserved, so that a
  coarse mesh reproduces the large scales of the direct simulation free of oscillations.
\item \emph{Why structure-preservation is necessary.} In the comparisons reported here, discretizations
  that fail to preserve the dissipation laws exhibit additional loss of the mean square potential and
  reach a different state at the same time, whereas the structure-preserving scheme reproduces selective decay and, in our
  experiments, approaches the Arnold-type bound.
\end{itemize}

The paper is organized as follows. \Cref{sec:model} writes the MHD-$\alpha$ model in the rotational
form used by the discretization and states its dissipation laws. \Cref{sec:turbulence-scheme}
presents the discretization and proves its structure-preserving properties, and
\Cref{sec:solver-development} describes the solver.
\Cref{sec:turbulence-numerical} verifies convergence and conservation.
\Cref{sec:selective-decay} studies selective decay and long-time relaxation,
\Cref{sec:spectra} compares the energy spectra with a direct simulation,
\Cref{sec:island} considers magnetic island coalescence, and
\Cref{sec:turbulence-conclusion} concludes.

\section{The MHD-$\alpha$ model and conserved quantities}\label{sec:model}
With the pressure $P=p+\ov{\bm u}\cdot\bm u-\frac S2|\bm B|^2$, the MHD-$\alpha$ model \eqref{eqn:mhd-alpha} reads, in rotational form,
\begin{subequations}\label{eqn:2d-mhd}
\begin{align}
&\partial_t \bm u - \ov{\bm u}\times w + \nu\curl w + \nabla P - S\, j\times\bm B = \bm f,\\
&\partial_t A + \ov{\bm u}\cdot\nabla A = -\eta j,\label{eqn:2d-mhd-A}\\
&\bm B = \curl A,\quad j = \rot\bm B,\quad w = \rot\bm u,\label{eqn:2d-mhd-B}\\
&\ov{\bm u} + \alpha^{2}\curl(\rot\ov{\bm u}) = \bm u,\quad \nabla\cdot\bm u = \nabla\cdot\ov{\bm u}=\nabla\cdot \bm B = 0.
\end{align}
\end{subequations}
For $\bm f=0$, the energy $\mathcal E^\alpha$ in \eqref{eqn:E-alpha}, the cross helicity $\mathcal H_c$ and the mean square potential $\mathcal A$ in \eqref{eqn:mhd-invariants} satisfy the dissipation laws
\begin{subequations}\label{eqn:2d-laws}
\begin{align}
\frac{d}{dt}\mathcal{E}^{\alpha} &= -\nu\int_{\Omega} w\,\ov{w}\dx - S\eta\int_{\Omega} |j|^{2}\dx,\\
\frac{d}{dt}\mathcal{H}_c &= -(\nu+\eta)\int_{\Omega} j\, w \dx,\\
\frac{d}{dt}\mathcal{A} &= -2\eta\int_{\Omega} A\, j \dx,
\end{align}
\end{subequations}
where $\ov{w} = \rot\ov{\bm u}$. For brevity we refer to the three quantities as invariants, even in the nonideal case.

\section{Structure-preserving finite element discretization}\label{sec:turbulence-scheme}

\subsection{Notations and finite element spaces}\label{ssec:spaces}
Let $\Omega=\mathbb T^2$ be the periodic square and let $\mathcal T_h$ be a shape-regular, quasi-uniform triangulation of $\mathbb T^2$ with mesh size $h$ (the cell diameter). We use $\|\cdot\|$ and $(\cdot,\cdot)$ for the $L^2(\Omega)$ norm and inner product, allowing $L^2$ to denote both scalar and vector-valued spaces. For $k\ge1$ let $V_h$, $\bm F_h$ and $Q_h$ be the Lagrange space of degree $k$, the Raviart--Thomas space \cite{raviart2006mixed} of degree $k$, and the mean-zero discontinuous space of degree $k-1$ on $\mathcal T_h$ respectively. These spaces are periodic by construction,
\begin{equation}
V_h\subset H^{1}(\mathbb T^2),\qquad
\bm F_h\subset \bm H(\div,\mathbb T^2),\qquad
Q_h\subset L^{2}_{0}(\mathbb T^2),
\end{equation}
and they form the finite element de\,Rham complex
\begin{equation}\label{eq:dRcomplex-2d}
\mathbb R\hookrightarrow V_h\xrightarrow{\,\curl\,}
                       \bm F_h\xrightarrow{\,\div\,}
                       Q_h\to 0,
\end{equation}
which is exact up to the harmonic forms of the torus:
\begin{equation}\label{eq:exact-2d}
\ker\curl|_{V_h}=\mathbb R,\qquad
\ker\div|_{\bm F_h}=\curl V_h\oplus\mathbb R^{2},\qquad
\div\bm F_h=Q_h,
\end{equation}
where $\mathbb R^{2}$ denotes the constant vector fields. We denote the discretely divergence-free subspace by $\bm Z_h:=\ker\div|_{\bm F_h}$. Since $\div\bm F_h\subseteq Q_h$, a divergence constraint imposed weakly against $Q_h$ holds pointwise. The discrete rotation $\rot_h:\bm F_h\to V_h$ is the $L^2$-adjoint of $\curl$,
\begin{equation}\label{eq:ibp-2d}
(\rot_h\bm v,\xi)=(\bm v,\curl\xi)
\qquad\forall\,\xi\in V_h,\ \bm v\in \bm F_h.
\end{equation}
Constant fields lie in the kernel of $\rot_h$, since $(\bm c,\curl\xi)=0$ for all $\xi\in V_h$, and taking $\xi=1$ in \eqref{eq:ibp-2d} shows that $\rot_h\bm v$ has zero mean for every $\bm v\in\bm F_h$. 


\begin{lemma}\label{lem:estimates}
There exist constants $c_p>0$, depending only on $\Omega$, and $C_{\mathrm{inv}}>0$, independent of $h$ and depending only on the shape regularity of $\mathcal T_h$ and on $k$, such that
\begin{align}
\|\bm v\|&\le c_p\|\rot_h\bm v\|\quad\forall\,\bm v\in\curl V_h, \label{eq:poincare-2d}\\
\|\curl\xi\|&\le C_{\mathrm{inv}}h^{-1}\|\xi\|\quad\forall\,\xi\in V_h.\label{eq:inv-curl-2d}
\end{align}
\end{lemma}
\begin{proof}
For \eqref{eq:poincare-2d}, write $\bm v=\curl\xi$ with $\xi\in V_h$. Since $\curl$ annihilates constants, we may subtract the mean of $\xi$ and assume that $\xi$ has zero mean, and this does not change $\bm v$. Then $\|\bm v\|^2=(\xi,\rot_h\bm v)\le c_p\|\nabla\xi\|\,\|\rot_h\bm v\|=c_p\|\bm v\|\,\|\rot_h\bm v\|$ by \eqref{eq:ibp-2d}, the Poincar\'e inequality for mean-zero $\xi$ and $|\curl\xi|=|\nabla\xi|$ pointwise. By \eqref{eq:exact-2d}, $\curl V_h$ consists of the mean-zero fields in $\bm Z_h$. The remaining two estimates are the standard inverse inequality on quasi-uniform meshes \cite[Thm.~4.5.11]{brenner2008mathematical}; for the last estimate, take $\xi=\rot_h\bm v$ in \eqref{eq:ibp-2d}.
\end{proof}


\subsection{Discrete filters and \texorpdfstring{$\alpha$}{alpha}-norms}\label{ssec:filters}
The continuous filter $\ov{\bm u}+\alpha^{2}\curl\rot\ov{\bm u}=\bm u$ is discretized in mixed form, with the auxiliary scalar $\ov w=\rot\ov{\bm u}$. The \emph{velocity filter} reads: given $\bm u\in\bm F_h$, find $(\ov{\bm u},\ov w)\in\bm F_h\times V_h$ such that
\begin{subequations}\label{eq:filter-2d}
\begin{align}
(\ov{\bm u},\bm\psi)+\alpha^{2}(\curl\ov w,\bm\psi)&=(\bm u,\bm\psi),\quad &&\forall\,\bm\psi\in\bm F_h,\label{eq:filter-2d-a}\\
(\ov w,\xi)-(\ov{\bm u},\curl\xi)&=0,\quad &&\forall\,\xi\in V_h.\label{eq:filter-2d-b}
\end{align}
\end{subequations}
By~\eqref{eq:ibp-2d}, equation~\eqref{eq:filter-2d-b} states $\ov w=\rot_h\ov{\bm u}$. The analysis of the viscous term also requires the \emph{vorticity filter} on $V_h$: given $\omega\in V_h$, find $\ov\omega\in V_h$ such that
\begin{equation}\label{eq:scalar-filter}
(\ov\omega,\varsigma)+\alpha^{2}(\curl\ov\omega,\curl\varsigma)
=(\omega,\varsigma)\qquad\forall\,\varsigma\in V_h.
\end{equation}

We denote by $\mathcal F_{\alpha}:\bm F_h\to\bm F_h$ and
$\mathcal R_{\alpha}:\bm F_h\to V_h$ the solution operators of the velocity filter,
$\bm u\mapsto\ov{\bm u}$ and $\bm u\mapsto\ov w$, and by
$\mathcal G_{\alpha}:V_h\to V_h$ the solution operator of the vorticity filter,
$\omega\mapsto\ov\omega$. All three are linear.

\begin{proposition}[Discrete filters]\label{prop:filters}
For every $\alpha>0$ the following hold.
\begin{enumerate}[label=(\roman*)]
\item Both filters are well posed, and
\begin{equation}\label{eq:filter-identity}
(\bm u,\ov{\bm u})=\|\ov{\bm u}\|^{2}+\alpha^{2}\|\ov w\|^{2},
\qquad
(\omega,\ov\omega)=\|\ov\omega\|^{2}+\alpha^{2}\|\curl\ov\omega\|^{2}.
\end{equation}
\item $\mathcal F_\alpha$ and $\mathcal G_\alpha$ are symmetric and contractive:
\begin{equation}\label{eq:filter-contractive}
\|\mathcal F_{\alpha}\bm u\|\le\|\bm u\|,
\qquad
\|\mathcal G_{\alpha}\omega\|\le\|\omega\|.
\end{equation}
\item The filters commute with $\rot_h$: $\mathcal R_\alpha\bm u=\mathcal G_\alpha(\rot_h\bm u)$, that is,
\begin{equation}\label{eq:filter-commute}
\rot_h\mathcal F_\alpha\bm u=\mathcal G_\alpha\rot_h\bm u\qquad\forall\,\bm u\in\bm F_h.
\end{equation}
\end{enumerate}
\end{proposition}

\begin{proof}
(i) The systems are finite-dimensional, so it suffices to show uniqueness. For the velocity filter, take $\bm\psi=\ov{\bm u}$ in \eqref{eq:filter-2d-a} and $\xi=\alpha^{2}\ov w$ in \eqref{eq:filter-2d-b}; adding, the cross terms cancel and $(\bm u,\ov{\bm u})=\|\ov{\bm u}\|^{2}+\alpha^{2}\|\ov w\|^{2}$, which vanishes when $\bm u=\bm0$. The vorticity filter is treated by taking $\varsigma=\ov\omega$ in \eqref{eq:scalar-filter}.

(ii) For $\bm u_1,\bm u_2\in\bm F_h$, taking $\bm\psi=\ov{\bm u}_2$ in \eqref{eq:filter-2d-a} for $\bm u_1$ and using \eqref{eq:filter-2d-b} for $\bm u_2$ gives $(\bm u_1,\ov{\bm u}_2)=(\ov{\bm u}_1,\ov{\bm u}_2)+\alpha^{2}(\ov w_1,\ov w_2)$, which is symmetric; $\mathcal G_\alpha$ is symmetric because \eqref{eq:scalar-filter} is. By \eqref{eq:filter-identity}, $\|\ov{\bm u}\|^{2}\le(\bm u,\ov{\bm u})\le\|\bm u\|\|\ov{\bm u}\|$, and likewise for $\ov\omega$.

(iii) For $\varsigma\in V_h$ take $\bm\psi=\curl\varsigma\in\bm F_h$ in \eqref{eq:filter-2d-a}. By \eqref{eq:ibp-2d} and $\ov w=\rot_h\ov{\bm u}$,
$(\ov w,\varsigma)+\alpha^{2}(\curl\ov w,\curl\varsigma)=(\rot_h\bm u,\varsigma)$, which is \eqref{eq:scalar-filter} with $\omega=\rot_h\bm u$.
\end{proof}

The discrete energy norm on $\bm F_h$ and dissipation norm on $V_h$ are
\begin{align}\label{eq:filter-norms}
&\|\bm v\|_{E,\alpha}^{2}:=(\bm v,\mathcal F_{\alpha}\bm v)=\|\ov{\bm v}\|^{2}+\alpha^{2}\|\mathcal R_\alpha\bm v\|^{2},
\qquad\\
&\|\omega\|_{D,\alpha}^{2}:=(\omega,\mathcal G_{\alpha}\omega)=\|\ov\omega\|^{2}+\alpha^{2}\|\curl\ov\omega\|^{2}.
\end{align}
By \Cref{prop:filters} these are induced by symmetric bilinear forms, and they are positive definite: $\|\bm v\|_{E,\alpha}=0$ forces $\ov{\bm v}=\bm0$, $\mathcal R_\alpha\bm v=0$, and then \eqref{eq:filter-2d-a} gives $\bm v=\bm 0$. The argument for $\|\cdot\|_{D,\alpha}$ is identical. 

\begin{theorem}[Norm equivalence]\label{thm:norm-equiv}
Let $\kappa>0$ be fixed and assume $\alpha\le\kappa h$. Then
\begin{equation}\label{eq:norm-equiv-2d}
c_E\|\bm v\|\le\|\bm v\|_{E,\alpha}\le\|\bm v\|\qquad\forall\,\bm v\in\bm F_h,
\end{equation}
where $c_E:=\bigl(\sqrt{2}\max\{1,C_{\mathrm{inv}}\kappa\}\bigr)^{-1}$, and $C_{\mathrm{inv}}$ is the constant of \Cref{lem:estimates}. In particular $c_E$ is independent of $h$ and $\alpha$.
\end{theorem}

\begin{proof}
Let $\bm v\in\bm F_h$ and write $\ov{\bm v}=\mathcal F_\alpha\bm v$, $\ov w=\mathcal R_\alpha\bm v$. For the upper bound, the Cauchy--Schwarz inequality and \eqref{eq:filter-contractive} give
\begin{equation}
\|\bm v\|_{E,\alpha}^{2}=(\bm v,\ov{\bm v})\le\|\bm v\|\,\|\ov{\bm v}\|\le\|\bm v\|^{2}.
\end{equation}
For the lower bound, take $\bm\psi=\bm v$ in \eqref{eq:filter-2d-a}. By the Cauchy--Schwarz inequality, the inverse estimate \eqref{eq:inv-curl-2d} and $\alpha^{2}h^{-1}\le\kappa\alpha$,
\begin{equation}
\begin{aligned}
\|\bm v\|^{2}&=(\ov{\bm v},\bm v)+\alpha^{2}(\curl\ov w,\bm v)\\
&\le\bigl(\|\ov{\bm v}\|+\alpha^{2}C_{\mathrm{inv}}h^{-1}\|\ov w\|\bigr)\|\bm v\|\\
&\le\bigl(\|\ov{\bm v}\|+C_{\mathrm{inv}}\kappa\,\alpha\|\ov w\|\bigr)\|\bm v\|.
\end{aligned}
\end{equation}
Dividing by $\|\bm v\|$, factoring out $\max\{1,C_{\mathrm{inv}}\kappa\}$ and using $a+b\le\sqrt2\,(a^{2}+b^{2})^{1/2}$,
\begin{equation}
\begin{aligned}
\|\bm v\|&\le\max\{1,C_{\mathrm{inv}}\kappa\}\bigl(\|\ov{\bm v}\|+\alpha\|\ov w\|\bigr)\\
&\le\sqrt2\max\{1,C_{\mathrm{inv}}\kappa\}\bigl(\|\ov{\bm v}\|^{2}+\alpha^{2}\|\ov w\|^{2}\bigr)^{1/2}\\
&=c_E^{-1}\|\bm v\|_{E,\alpha}.
\end{aligned}
\end{equation}
\end{proof}

\subsection{The fully discrete scheme}\label{ssec:scheme}
The initial data are taken as $\bm u_h^0=\nabla^{\perp}I_h\varphi\in\bm Z_h$ and  $A_h^0=-I_h\psi+\lvert\Omega\rvert^{-1}\int_\Omega I_h\psi\,\dx$ so that $\int_\Omega A_h^0\,\dx=0$, where $\nabla^{\perp}f:=(-\partial_y f,\ \partial_x f)=-\curl f$, $I_h$ is the Lagrange interpolant into $V_h$, and $\varphi,\psi$ are the stream functions of the velocity and the magnetic field, so that $\bm u_h^0$ and $\bm B_h^0=\curl A_h^0=\nabla^{\perp}I_h\psi$ are discretely divergence-free by construction. The filter width $\alpha$ is constant on the mesh. 

We present a structure-preserving discretization in space and time for \eqref{eqn:2d-mhd}. We use the implicit midpoint method as the temporal scheme. For the fluid variables and the potential we use the averages
\begin{equation}
    \bm u_h^{n+1/2} := \tfrac{1}{2}(\bm u_h^{n+1} + \bm u_h^n),
    \quad
    \ov{\bm u}_h^{n+1/2} := \tfrac{1}{2}(\ov{\bm u}_h^{n+1} + \ov{\bm u}_h^n),
    \quad
    A_h^{n+1/2} := \tfrac{1}{2}(A_h^{n+1}+A_h^n).
\end{equation}
The variables $w_h^{n+1/2}$, $j_h^{n+1/2}$ and $P_h^{n+1/2}$ are treated as independent unknowns at the intermediate time level. The filter constraints for $\ov{\bm u}_h^{n+1}$ and $\ov w_h^{n+1}$ are imposed at the endpoint $t^{n+1}$, so that the energy pairing $(\bm u_h^{n+1}, \ov{\bm u}_h^{n+1})$ is evaluated at a single time level.

The formulation is chosen as follows. Conservation of the mean square potential relies on $\div\ov{\bm u}_h=0$ exactly, since exact divergence-free ensures that $(\ov{\bm u}_h\cdot\nabla A_h,A_h)=0$, so we discretize $\ov{\bm u}_h$ in the $\bm H(\div)$-conforming space $\bm F_h$. However, $\bm F_h$ is not regular enough for the filter $\ov{\bm u}+\alpha^{2}\curl\rot\ov{\bm u}=\bm u$, so we write it in mixed form with the auxiliary variable $\ov w=\rot\ov{\bm u}$. Since $\bm u_h=\ov{\bm u}_h+\alpha^{2}\curl\ov w_h$ and $\div\curl=0$, $\div\bm u_h=\div\ov{\bm u}_h$, so $\bm u_h$ must also be
exactly divergence free, and we use the pair $\bm F_h\times Q_h$ for the velocity and
pressure, for which $\div\bm F_h=Q_h$ makes the discretely divergence-free condition hold pointwise. The viscous term is written in
mixed form with $w=\rot\bm u$ thanks to the periodic boundary conditions. For the energy and the cross helicity we follow \cite{huHelicityconservativeFiniteElement2021}, with one modification: since the momentum equation is tested with $\ov{\bm u}_h$ instead of $\bm u_h$, the kinetic energy is $\frac{1}{2}(\bm u_h,\ov{\bm u}_h)$, which requires $\ov{\bm u}_h^{n}=\mathcal F_\alpha\bm u_h^{n}$ at both time levels, so $\ov{\bm u}_h^{n+1}$ and $\ov w_h^{n+1}$ are included as independent unknowns of \Cref{alg:sp-discretisation-mhd-alpha}, with the filter equations imposed at $t^{n+1}$.

\begin{algorithm}[H]
\caption{Fully discrete scheme for the 2D incompressible MHD-$\alpha$ model}\label{alg:sp-discretisation-mhd-alpha}
Given $(\bm u_h^0,\, A_h^0)\in\bm Z_h\times V_h$, set $\ov{\bm u}_h^0:=\mathcal F_\alpha\bm u_h^0$.
For $n=0,1,\ldots,N-1$, find
\begin{equation}
\begin{gathered}
(\bm u_h^{n+1},\, \ov{\bm u}_h^{n+1},\,
 P_h^{n+1/2},\, w_h^{n+1/2},\, \ov{w}_h^{n+1},\, A_h^{n+1},\, j_h^{n+1/2})\\
\in [\bm F_h]^2 \times Q_h \times [V_h]^4
\end{gathered}
\end{equation}
such that for all test functions
$(\bm v_h,\, \bm\chi_h,\, q_h,\, z_h,\, \xi_h,\, e_h,\, k_h)$
in the same space:
\begin{subequations}\label{eq:sp-scheme-2d}
\begin{empheq}[left=\empheqlbrace]{align}
&\frac{1}{\Delta t}(\bm u_h^{n+1} - \bm u_h^n,\, \bm v_h)
 - (\ov{\bm u}_h^{n+1/2} \times w_h^{n+1/2},\, \bm v_h)
 + \nu(\curl w_h^{n+1/2},\, \bm v_h) \notag\\
&\qquad - (P_h^{n+1/2},\, \div\bm v_h)
 - S(j_h^{n+1/2}\times\bm B_h^{n+1/2},\, \bm v_h) \notag\\
&\qquad = (\bm f(t^{n+1/2}),\, \bm v_h),\\
&(\ov{\bm u}_h^{n+1},\, \bm\chi_h)
 + \alpha^2(\curl \ov{w}_h^{n+1},\, \bm\chi_h)
 = (\bm u_h^{n+1},\, \bm\chi_h),\\
&(\div\bm u_h^{n+1},\, q_h) = 0,\\
&(w_h^{n+1/2},\, z_h) = (\bm u_h^{n+1/2},\, \curl z_h),\\
&(\ov{w}_h^{n+1},\, \xi_h) = (\ov{\bm u}_h^{n+1},\, \curl \xi_h),\\
&\frac{1}{\Delta t}(A_h^{n+1} - A_h^n,\, e_h)
 + \eta(j_h^{n+1/2},\, e_h) \notag\\
&\qquad - (\ov{\bm u}_h^{n+1/2}\times \bm B_h^{n+1/2},\, e_h) = 0, \label{eq:sp-scheme-2d-A}\\
&(j_h^{n+1/2},\, k_h) = (\bm B_h^{n+1/2},\, \curl k_h).\label{eq:sp-scheme-2d-j}
\end{empheq}
\end{subequations}
\end{algorithm}

Throughout, $\bm B_h^{n} := \curl A_h^{n}$, so that
$\bm B_h^{n+1/2} := \curl A_h^{n+1/2} = \tfrac12(\bm B_h^{n+1}+\bm B_h^{n})$
and $\div \bm B_h^n = 0$. The second and fifth equations of \Cref{alg:sp-discretisation-mhd-alpha} are the velocity
filter \eqref{eq:filter-2d} at $t^{n+1}$, so $\ov{\bm u}_h^{n}=\mathcal F_\alpha\bm u_h^{n}$
and, by linearity, $\ov{\bm u}_h^{n+1/2}=\mathcal F_\alpha\bm u_h^{n+1/2}$. The fourth and
last equations state $w_h^{n+1/2}=\rot_h\bm u_h^{n+1/2}$ and
$j_h^{n+1/2}=\rot_h\bm B_h^{n+1/2}$.

\subsection{Structure-preserving properties}\label{ssec:sp-properties}

For convenience, we write $D_t\bm u_h:=(\bm u_h^{n+1}-\bm u_h^n)/\Delta t$, and define $D_tA_h$ and $D_t\bm B_h$ analogously.

\begin{theorem}\label{thm:2d-sp}

Every solution of the discrete system \eqref{eq:sp-scheme-2d} satisfies the following for all $n$.
\begin{enumerate}[label=(\alph*)]
\item \textit{discretely divergence-free property.}
\begin{equation}\label{eqn:div-free-disc-2d}
\div \bm u_h^n = \div \ov{\bm u}_h^n = \div \bm B_h^n = 0,
\end{equation}
provided the initial conditions are discretely divergence-free.
\item \textit{discrete energy dissipation law}, with $\mathcal{E}_h^{\alpha, n} := \frac12\bigl(\|\bm u_h^n\|_{E,\alpha}^{2} + S\|\bm B_h^n\|^{2}\bigr)$,
\begin{equation}\label{eqn:energy-disc-2d}
    \mathcal{E}_h^{\alpha, n+1} - \mathcal{E}_h^{\alpha, n}
        = -\Delta t\,\nu\,\|w_h^{n+1/2}\|^2_{D,\alpha}
          -\Delta t\,S\eta\,\|j_h^{n+1/2}\|^{2}
          +\Delta t\,(\bm f^{n+1/2}, \ov{\bm u}_h^{n+1/2}).
\end{equation}
In the ideal limit $\nu = \eta = 0$ with $\bm f = \bm 0$, $\mathcal{E}_h^{\alpha}$ is exactly conserved. If $\nu>0$ and $\int_\Omega\bm f(t^{n+1/2})\,\dx=\bm0$ for every $n$, then for every $N\ge 1$,
\begin{equation}\label{eq:stability-2d}
\mathcal{E}_h^{\alpha, N}
+\frac{\nu\Delta t}{2}\sum_{n=0}^{N-1}\|w^{n+1/2}_h\|_{D,\alpha}^2
 +S\eta\,\Delta t\sum_{n=0}^{N-1}\|j^{n+1/2}_h\|^{2}
\leq \mathcal{E}_h^{\alpha, 0}
   +\frac{c_p^{2}\Delta t}{2\nu}\sum_{n=0}^{N-1}\|\bm f^{n+1/2}\|^2,
\end{equation}
so the scheme is unconditionally stable.
\item \textit{discrete cross helicity dissipation law}, with $\bm f=\bm 0$ and $\mathcal{H}_{c,h}^{n} := (\bm u_h^n, \bm B_h^n)$,
\begin{equation}\label{eqn:cross-helicity-disc-2d}
    \mathcal{H}_{c,h}^{n+1} - \mathcal{H}_{c,h}^{n}
        = -\Delta t\,(\nu + \eta)\,(j_h^{n+1/2}, w_h^{n+1/2}).
\end{equation}
\item \textit{discrete mean-square potential dissipation law}, with $\mathcal{A}_h^n := \|A_h^n\|^{2}$,
\begin{equation}\label{eqn:msp-disc-2d}
    \mathcal{A}_h^{n+1} - \mathcal{A}_h^n
        = -2\Delta t\,\eta\,\|\nabla A_h^{n+1/2}\|^{2}.
\end{equation}
\end{enumerate}
In the ideal limit, $\mathcal{H}_{c,h}$ is conserved for $\nu=\eta=0$ and $\mathcal{A}_h$ for $\eta=0$.
\end{theorem}

\begin{proof}
Throughout the proof we drop the superscript $n+1/2$ on midpoint quantities. By \Cref{alg:sp-discretisation-mhd-alpha} we have $\bm B_h=\curl A_h$, $w_h=\rot_h\bm u_h$, $j_h=\rot_h\bm B_h$ and $\ov{\bm u}_h=\mathcal F_\alpha\bm u_h$, so that \eqref{eq:ibp-2d} gives
\begin{equation}
(D_tA_h,j_h)=(D_t\bm B_h,\bm B_h),\qquad (D_tA_h,w_h)=(D_t\bm B_h,\bm u_h).
\end{equation}
We also use the pointwise identity $(\bm a\times\phi)\cdot\bm b=-\phi\,(\bm a\times\bm b)$ for $\bm a,\bm b\in\mathbb R^2$, $\phi\in\mathbb R$, which implies $(\bm a\times\phi)\cdot\bm a=0$ and
\begin{equation}
-(j_h\times\bm B_h,\ov{\bm u}_h)=(\ov{\bm u}_h\times\bm B_h,j_h),\qquad
-(\ov{\bm u}_h\times w_h,\bm B_h)=(\ov{\bm u}_h\times\bm B_h,w_h).
\end{equation}

The divergence-free properties of $\bm u_h$ and $\bm B_h$ are due to the properties of the de Rham complex. Since \eqref{eq:filter-2d-a} holds for all $\bm\psi\in\bm F_h$ and all terms lie in $\bm F_h$, $\ov{\bm u}_h^{n}=\bm u_h^{n}-\alpha^{2}\curl\ov w_h^{n}$ is a sum of two divergence-free fields.

For the energy law, we test the momentum equation with $\ov{\bm u}_h$ and the potential equation with $e_h=Sj_h$. The nonlinear term drops since $(\ov{\bm u}_h\times w_h)\cdot\ov{\bm u}_h=0$, and since $\div\ov{\bm u}_h=0$ the pressure term vanishes, so that
\begin{align}
(D_t\bm u_h,\ov{\bm u}_h)+\nu(\curl w_h,\ov{\bm u}_h)-S(j_h\times\bm B_h,\ov{\bm u}_h)&=(\bm f,\ov{\bm u}_h),\\
S(D_t\bm B_h,\bm B_h)+S\eta\|j_h\|^2-S(\ov{\bm u}_h\times\bm B_h,j_h)&=0.
\end{align}
Adding the two equations, the coupling terms cancel by the pointwise identity, and it remains to identify the other terms. Since $\mathcal F_\alpha$ is symmetric by \Cref{prop:filters},
\begin{equation}
(D_t\bm u_h,\ov{\bm u}_h)=\frac{\|\bm u_h^{n+1}\|_{E,\alpha}^2-\|\bm u_h^{n}\|_{E,\alpha}^2}{2\Delta t},
\end{equation}
and by \eqref{eq:ibp-2d} and \eqref{eq:filter-commute},
\begin{equation}
\nu(\curl w_h,\ov{\bm u}_h)=\nu(w_h,\rot_h\mathcal F_\alpha\bm u_h)=\nu(w_h,\mathcal G_\alpha w_h)=\nu\|w_h\|_{D,\alpha}^2.
\end{equation}
Multiplying by $\Delta t$ gives \eqref{eqn:energy-disc-2d}.

To obtain the stability bound, we estimate the forcing term. Let $\bm c$ be the mean of $\ov{\bm u}_h\in\bm Z_h$, so that $\ov{\bm u}_h-\bm c\in\curl V_h$. Since $\bm f$ has zero mean and $\rot_h\bm c=0$, Young's inequality, \eqref{eq:poincare-2d} and $\rot_h\ov{\bm u}_h=\mathcal G_\alpha w_h$ give
\begin{equation}
|(\bm f,\ov{\bm u}_h)|=|(\bm f,\ov{\bm u}_h-\bm c)|
\le\tfrac12c_p^{2}\nu^{-1}\|\bm f\|^2+\tfrac12\nu c_p^{-2}\|\ov{\bm u}_h-\bm c\|^2
\le\tfrac12c_p^{2}\nu^{-1}\|\bm f\|^2+\tfrac12\nu\|w_h\|_{D,\alpha}^2.
\end{equation}
Inserting this into \eqref{eqn:energy-disc-2d} and summing over $n$ gives \eqref{eq:stability-2d}.

Next, for the cross helicity, we test the momentum equation with $\bm B_h$ and the potential equation with $e_h=w_h$. Since $\div\bm B_h=0$ the pressure term vanishes, $(j_h\times\bm B_h,\bm B_h)=0$, $\nu(\curl w_h,\bm B_h)=\nu(w_h,j_h)$ and $(D_tA_h,w_h)=(D_t\bm B_h,\bm u_h)$, so that
\begin{align}
(D_t\bm u_h,\bm B_h)-(\ov{\bm u}_h\times w_h,\bm B_h)+\nu(w_h,j_h)&=0,\\
(D_t\bm B_h,\bm u_h)+\eta(j_h,w_h)-(\ov{\bm u}_h\times\bm B_h,w_h)&=0.
\end{align}
Adding the two equations, the nonlinear terms cancel by the pointwise identity, and since
$(D_t\bm u_h,\bm B_h^{n+1/2})+(D_t\bm B_h,\bm u_h^{n+1/2})=\Delta t^{-1}(\mathcal{H}_{c,h}^{n+1}-\mathcal{H}_{c,h}^{n})$, we obtain \eqref{eqn:cross-helicity-disc-2d}.

Finally, for the mean-square potential, we test the potential equation with $e_h=A_h^{n+1/2}$. Since $\div\ov{\bm u}_h=0$ and $\Omega$ is periodic,
\begin{equation}
-(\ov{\bm u}_h\times\bm B_h,A_h)=(\ov{\bm u}_h\cdot\nabla A_h,A_h)=\int_\Omega\nabla\cdot\Bigl(\ov{\bm u}_h\frac{A_h^2}{2}\Bigr)\dx=0,
\end{equation}
so the equation reduces to $\frac{1}{2\Delta t}(\|A_h^{n+1}\|^2-\|A_h^n\|^2)+\eta(j_h,A_h)=0$. Since $(j_h,A_h)=(\bm B_h,\curl A_h)=\|\nabla A_h\|^2$ by \eqref{eq:ibp-2d}, this is \eqref{eqn:msp-disc-2d}.
\end{proof}


\begin{remark}[$L^2$ stability]\label{rem:L2-stability}
Under the assumptions of \Cref{thm:2d-sp}(b) and \Cref{thm:norm-equiv}, \eqref{eq:stability-2d} and \eqref{eq:norm-equiv-2d} give
$c_E^{2}\|\bm u_h^N\|^{2}+S\|\bm B_h^N\|^{2}
\le \|\bm u^{0}_h\|^{2}+S\|\bm B^{0}_h\|^{2}
+\Delta t\,c_p^{2}\nu^{-1}\sum_{n=0}^{N-1}\|\bm f^{n+1/2}\|^2$,
with $c_E$ independent of $h$, $\alpha$ and $\Delta t$.
\end{remark}

\section{Nonlinear solver and block preconditioning}
\label{sec:solver-development}
Turbulence simulations typically require fine meshes and long time computation. Therefore, we develop block preconditioners for the fully discrete scheme of \Cref{alg:sp-discretisation-mhd-alpha}. At each time step the nonlinear problem of \Cref{alg:sp-discretisation-mhd-alpha} is solved by Newton's method. We group the pressure separately from the other fields and apply an augmented Lagrangian/grad-div stabilization term \cite{benziAugmentedLagrangianBased2006,olshanskii2004grad}, so that the Jacobian system has the form
\begin{equation}
 \mathcal J_\gamma=
 \begin{bmatrix}
 \mathcal K_\gamma & \mathcal D^{\mathsf T}\\
 \mathcal D & 0
 \end{bmatrix},
 \qquad
 \mathcal K_\gamma=\mathcal K+\gamma \mathcal D^{\mathsf T}\mathcal M_p^{-1}\mathcal D,
 \label{eq:solver-jacobian}
\end{equation}
where $\mathcal K$ couples the six remaining fields, $\mathcal D$ is the discrete divergence and $\mathcal M_p$ is the pressure mass matrix. We take $\gamma=10^2$ in all experiments. We apply FGMRES with the upper triangular preconditioner
\begin{equation}
 \mathcal P_U=
 \begin{bmatrix}
 \widetilde{\mathcal K}_\gamma & \mathcal D^{\mathsf T}\\
 0 & \widetilde{\mathcal S}_\gamma
 \end{bmatrix},
 \qquad
 \widetilde{\mathcal S}_\gamma=-\gamma^{-1}\mathcal M_p,
 \label{eq:upper-schur}
\end{equation}
where $\widetilde{\mathcal S}_\gamma$ approximates the Schur complement $-\mathcal D\mathcal K_\gamma^{-1}\mathcal D^{\mathsf T}$ and is inverted by one Jacobi step. Since the experiments using this preconditioner use piecewise-constant pressure, $\mathcal M_p$ is diagonal and this step is exact.
For the inversion of the block $\widetilde{\mathcal K}_\gamma$, we apply the monolithic multigrid in \cite{laakmannAugmentedLagrangianPreconditioner2022}. The block $\widetilde{\mathcal K}_\gamma$ is solved inexactly with a fixed amount of work: we perform exactly two FGMRES iterations, preconditioned by one multigrid F-cycle, which starts with an exact solve on the coarsest mesh and performs one V-cycle on each finer level \cite{laakmannAugmentedLagrangianPreconditioner2022}. 

\begin{figure}[htbp]
\centering
\begin{forest}
  for tree={
    grow'=0,
    child anchor=west,
    parent anchor=south,
    anchor=west,
    draw,
    font=\small,
    inner xsep=5pt, inner ysep=4pt,
    s sep=9pt,
    calign=first,
    edge path={
      \noexpand\path[draw, \forestoption{edge}]
        (!u.south west) +(9pt,0) |- (.child anchor)\forestoption{edge label};
    },
    before typesetting nodes={if n=1{insert before={[,phantom]}}{}},
    fit=band,
    before computing xy={l=20pt},
  }
  [{Newton iteration}
    [{Krylov solver (FGMRES)}
      [{Upper block-triangular preconditioner $\mathcal P_U$}
        [{Schur complement approximation $\widetilde{\mathcal S}_\gamma=-\gamma^{-1}\mathcal M_p$}
          [{Exact pressure mass matrix inverse}]
        ]
        [{Solver for the augmented block $\widetilde{\mathcal K}_\gamma$}
          [{Krylov solver (FGMRES)}
            [{Monolithic F-cycle multigrid}
              [{Coarse grid solver}
                [{LU factorization}]
              ]
              [{Relaxation}
                [{GMRES}
                  [{Additive star iteration}]
                ]
              ]
            ]
          ]
        ]
      ]
    ]
  ]
\end{forest}
\caption{Graphical outline of the nonlinear solver and the block preconditioner for one time step of \Cref{alg:sp-discretisation-mhd-alpha}.}
\label{fig:solver-outline}
\end{figure}
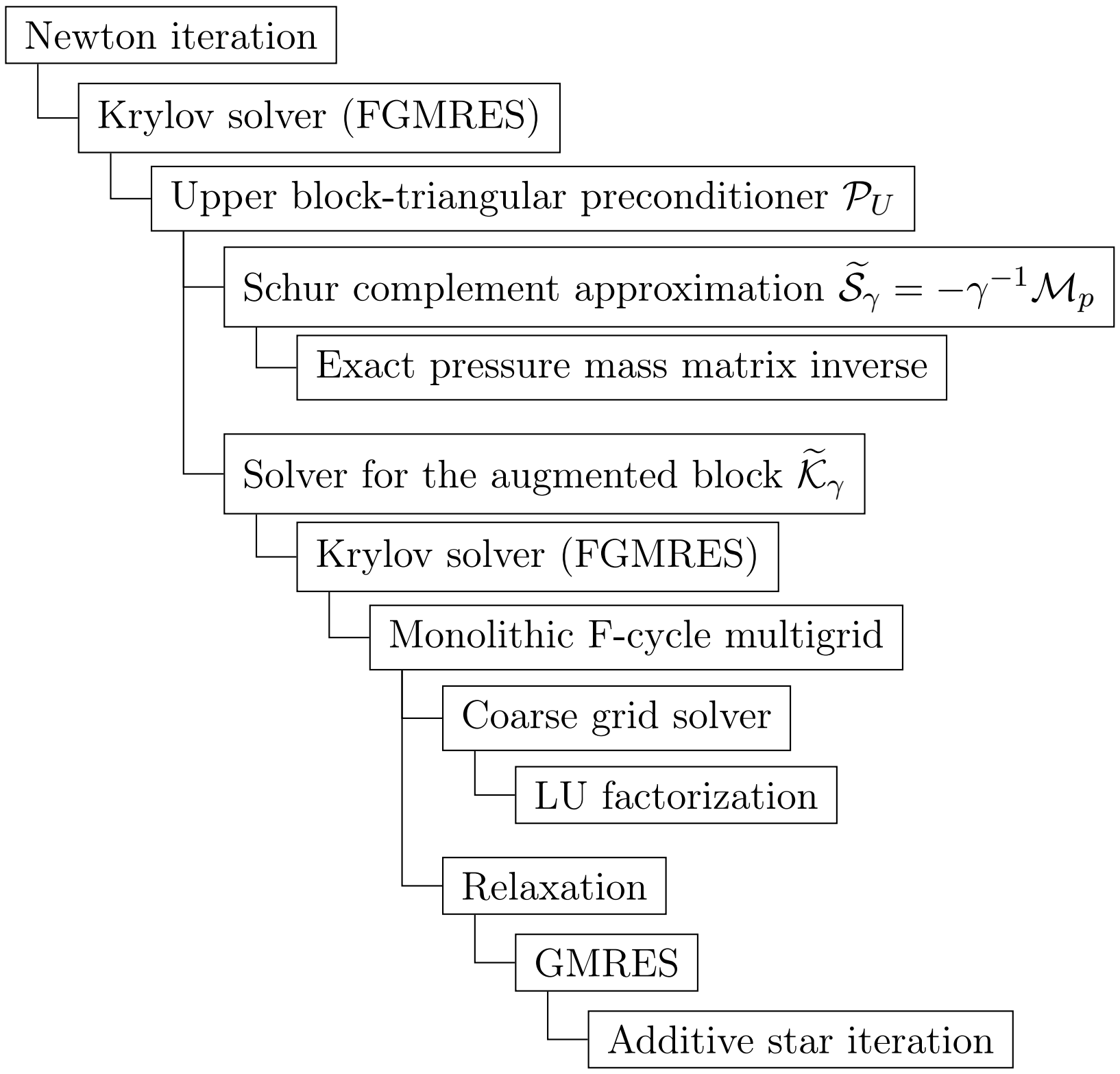

On the coarsest level we use a sparse LU factorization with MUMPS \cite{amestoy2001}.
On the finer levels the smoother is six GMRES iterations preconditioned by additive Schwarz over vertex-star patches, with the patch problems solved exactly by LU \cite{farrell2019cpcpatch}. We stop the nonlinear solver at a relative residual of $10^{-8}$, and the outer FGMRES iteration at a relative residual of $10^{-7}$ or an absolute residual of $10^{-10}$, measured in the Euclidean norm. The filter width $\alpha=c\,h_{\mathrm{fine}}$ is fixed on the finest mesh and transferred to the coarse levels. The performance of the solver is assessed in the numerical experiments below. The overall structure of the solver is summarized in \Cref{fig:solver-outline}.

\section{Numerical experiments}\label{sec:turbulence-numerical}
In this section, we present numerical results for the structure-preserving scheme of \Cref{sec:turbulence-scheme}. We first verify the scheme by a convergence test with manufactured solutions (\Cref{ssec:convergence}) and a conservation test on the Orszag--Tang vortex (\Cref{ssec:turbulence-conservation}). The scheme is then used to study selective decay, energy spectra and island coalescence in \Cref{sec:selective-decay,sec:spectra,sec:island}. The convergence and conservation tests, the Orszag--Tang run toward the Arnold-type bound (\Cref{ssec:ozt-arnold}) and the random-data experiment (\Cref{sec:random-relaxation}) solve the linear system in each Newton step by a sparse direct LU factorization (MUMPS \cite{amestoy2001}). All other experiments use the solver of \Cref{sec:solver-development}, and its robustness with respect to the physical and discretization parameters is reported in \Cref{tab:newton-performance,tab:reynolds-robustness}. 

\subsection{Convergence}\label{ssec:convergence}
\paragraph{Spatial convergence}
For the convergence tests only, we augment the momentum and magnetic-potential equations in \eqref{eq:sp-scheme-2d} by manufactured right-hand sides $(\bm f_{\rm ex},\bm v_h)$ and $(g_{A,{\rm ex}},e_h)$, respectively, obtained by substituting the exact fields into the continuous equations.
We test the spatial accuracy on the periodic square $\Omega = (0, 2\pi)^2$ with $T = 1$, $\Delta t = 0.1$ and physical parameters $\alpha = 0.1$, $\nu = \eta = S = 1$. We choose a stationary manufactured solution to focus only on the spatial discretization error.
The manufactured solution is
\begin{align}
\bm u_{\rm ex}(x,y) &= A_0\,(-\sin x \cos y,\; \cos x \sin y)^T, \\
A_{\rm ex}(x,y)     &= A_0\,\sin x \sin y, \\
P_{\rm ex}(x,y)     &= A_0\,\sin x \sin y,
\end{align}
with $A_0=0.1$, $\ov{\bm u}_{\rm ex} = \bm u_{\rm ex} / (1 + 2\alpha^2)$ and $\bm B = \curl A$ following analytically. Errors are measured in the $L^2$ norm for $\bm u, P, \bm B$ and in the $H^1$ norm for $A$. \Cref{tab:2d-spatial-convergence-k1,tab:2d-spatial-convergence-k2} report optimal $\mathcal{O}(h^k)$ rates for $k = 1, 2$.
\begin{table}[htbp]
    \centering
    \footnotesize
    \setlength{\tabcolsep}{3pt}
    \begin{tabular}{lcccccccc}
        \toprule
        $N$ & $\|e_u\|_{L^2}$ & Rate & $\|e_P\|_{L^2}$ & Rate & $\|e_B\|_{L^2}$ & Rate & $\|e_A\|_{H^1}$ & Rate \\
        \midrule
        $32$    & $4.367\times10^{-2}$ & --   & $2.054\times10^{-2}$ & --   & $4.349\times10^{-2}$ & --   & $4.364\times10^{-2}$ & --   \\
        $64$    & $2.182\times10^{-2}$ & 1.00 & $1.028\times10^{-2}$ & 1.00 & $2.179\times10^{-2}$ & 1.00 & $2.181\times10^{-2}$ & 1.00 \\
        $128$   & $1.091\times10^{-2}$ & 1.00 & $5.140\times10^{-3}$ & 1.00 & $1.090\times10^{-2}$ & 1.00 & $1.090\times10^{-2}$ & 1.00 \\
        \bottomrule
    \end{tabular}
     \caption{2D spatial convergence rate for $k=1$ on $N\times N$ meshes.}
    \label{tab:2d-spatial-convergence-k1}
\end{table}

\begin{table}[htbp]
    \centering
    \footnotesize
    \setlength{\tabcolsep}{3pt}
   
    \begin{tabular}{lcccccccc}
        \toprule
        $N$ & $\|e_u\|_{L^2}$ & Rate & $\|e_P\|_{L^2}$ & Rate & $\|e_B\|_{L^2}$ & Rate & $\|e_A\|_{H^1}$ & Rate \\
        \midrule
        $8$    & $2.579\times10^{-2}$ & --   & $1.300\times10^{-2}$ & --   & $2.575\times10^{-2}$ & --   & $2.588\times10^{-2}$ & --   \\
        $16$   & $6.684\times10^{-3}$ & 1.95 & $3.167\times10^{-3}$ & 2.04 & $6.679\times10^{-3}$ & 1.95 & $6.688\times10^{-3}$ & 1.95 \\
        $32$   & $1.686\times10^{-3}$ & 1.99 & $7.845\times10^{-4}$ & 2.01 & $1.685\times10^{-3}$ & 1.99 & $1.686\times10^{-3}$ & 1.99 \\
        \bottomrule
    \end{tabular}
     \caption{2D spatial convergence rate for $k=2$ on $N\times N$ meshes.}
    \label{tab:2d-spatial-convergence-k2}
\end{table}

\paragraph{Temporal convergence}
We test the temporal accuracy on the periodic square $\Omega = (0, 2\pi)^2$ with a fixed $96 \times 96$ mesh, $T = 1$,
$\alpha = 0.1$ and $\nu = \eta = S = 1$. We use polynomial degree $k=3$ in space, so that the spatial error is
negligible, and vary the time step. The manufactured solution is
\begin{align}
\bm u_{\rm ex} &= \curl(\sin x\sin y)\,\cos 2t, \quad\\
A_{\rm ex}     &= \bigl(\sin x\sin y + \tfrac12\cos x\bigr)\, e^{-t}, \quad\\
P_{\rm ex}     &= \cos x\cos y\,\cos t,
\end{align}
with $\ov{\bm u}_{\rm ex} = \bm u_{\rm ex}/(1 + 2\alpha^2)$ and $\bm B = \curl A$. The corresponding momentum and magnetic-potential source terms are
evaluated at $t^{n+1/2}$. The errors of $\bm u_h$, $\curl A_h$ and $A_h$ are measured at $t=T$, and the error of the
pressure, which is a midpoint quantity, at $t=T-\Delta t/2$. \Cref{tab:2d-temporal-convergence} shows the expected
$\mathcal{O}(\Delta t^2)$ rates of the implicit midpoint rule.
\begin{table}[htbp]
    \centering
    \footnotesize
    \setlength{\tabcolsep}{3pt}
    \begin{tabular}{lcccccccc}
        \toprule
        $\Delta t$ & $\|e_u\|_{H(\mathrm{div})}$ & Rate & $\|e_P\|_{L^2}$ & Rate & $\|e_B\|_{H(\mathrm{div})}$ & Rate & $\|e_A\|_{H^1}$ & Rate \\
        \midrule
        $1/4$  & $7.806\times10^{-3}$ & --   & $4.986\times10^{-3}$ & --   & $1.445\times10^{-2}$ & --   & $1.788\times10^{-2}$ & --   \\
        $1/8$    & $2.222\times10^{-3}$ & 1.81 & $2.458\times10^{-3}$ & 1.02 & $3.575\times10^{-3}$ & 2.01 & $4.424\times10^{-3}$ & 2.01 \\
        $1/16$   & $5.721\times10^{-4}$ & 1.96 & $7.829\times10^{-4}$ & 1.65 & $8.917\times10^{-4}$ & 2.00 & $1.103\times10^{-3}$ & 2.00 \\
        $1/32$   & $1.444\times10^{-4}$ & 1.99 & $2.174\times10^{-4}$ & 1.85 & $2.231\times10^{-4}$ & 2.00 & $2.759\times10^{-4}$ & 2.00 \\
        $1/64$   & $3.710\times10^{-5}$ & 1.96 & $5.725\times10^{-5}$ & 1.93 & $5.671\times10^{-5}$ & 1.98 & $6.973\times10^{-5}$ & 1.98 \\
        \bottomrule
    \end{tabular}
     \caption{2D temporal convergence test on the periodic square.}
    \label{tab:2d-temporal-convergence}
\end{table}

\subsection{Conservation test}\label{ssec:turbulence-conservation}
We now test the conservation properties of the scheme. The setting is the standard Orszag--Tang vortex on the periodic square $\Omega = (0,2\pi)^{2}$, in the ideal regime $\nu = \eta = 0$, so that the three quadratic invariants should be conserved exactly.
The initial magnetic and velocity fields are prescribed
in terms of scalar stream functions to ensure divergence-free conditions, as in \cite{biskamp1989dynamics}. The initial magnetic and velocity fields are defined by
\begin{align}
\mathbf{B}_0 &= \bigl(-\partial_y \psi, \; \partial_x \psi \bigr),\\
\mathbf{u}_0
&= \bigl(-\partial_y \varphi, \; \partial_x \varphi \bigr),
\end{align}
which satisfy $\nabla \cdot \mathbf{B}_0 = 0$ and $\nabla\cdot \bm u_0=0$. The scalar potentials $\varphi(x,y)$ and $\psi(x,y)$ are chosen as
\begin{align}
\varphi(x,y) &= \cos(x + 1.4) + \cos(y + 0.5), \\
\psi(x,y) &= \cos(2x + 2.3) + \cos(y + 4.1).
\end{align}

\begin{figure}[htbp]
    \centering
    \begin{subfigure}[b]{0.48\textwidth}
        \centering
        \includegraphics[width=\linewidth]{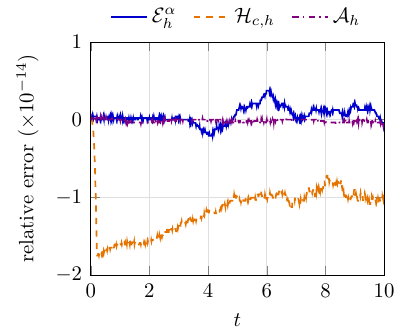}
        \caption{}
    \end{subfigure}
    \hfill
    \begin{subfigure}[b]{0.48\textwidth}
        \centering
        \includegraphics[width=\linewidth]{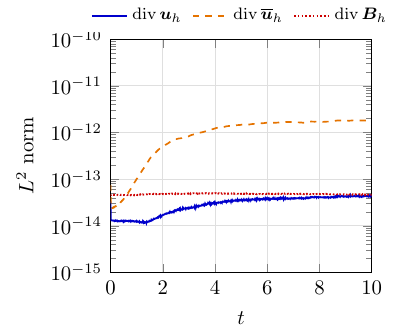}
        \caption{}
    \end{subfigure}
    \caption{Ideal Orszag--Tang test on the $64\times64$ mesh with $\alpha=h$ and $\Delta t=0.01$. (a) Relative errors $(Q_h(t)-Q_h(0))/Q_h(0)$ of the energy, cross helicity and mean square potential. (b) $L^2$ norms of $\div\ov{\bm u}_h$, $\div\bm u_h$ and $\div\bm B_h$.}
    \label{fig:sp2d-errors}
\end{figure}

We choose $\Delta t= 0.01$, $T = 10$ with $\nu = \eta = 0$, $S=1$, $\alpha = h$, where $h$ is the mesh size. The mesh is $64\times 64$ with the lowest-order element in space. The relative errors of the total energy, cross helicity and mean square potential are shown in \Cref{fig:sp2d-errors}. All quantities are preserved. 

\section{Selective decay}\label{sec:selective-decay}
In two-dimensional MHD, energy and cross helicity cascade forward to small scales, where they are dissipated by $\nu$ and $\eta$, while the mean square potential cascades inversely to large scales \cite{biskamp1989dynamics,pouquet1978two}. We now investigate how the scheme reproduces the resulting selective decay of the three invariants.
\subsection{Orszag--Tang: timescale separation and the Arnold-type bound}\label{ssec:ozt-arnold}
\begin{figure}[htpb]
    \centering
    \includegraphics[width=\textwidth]{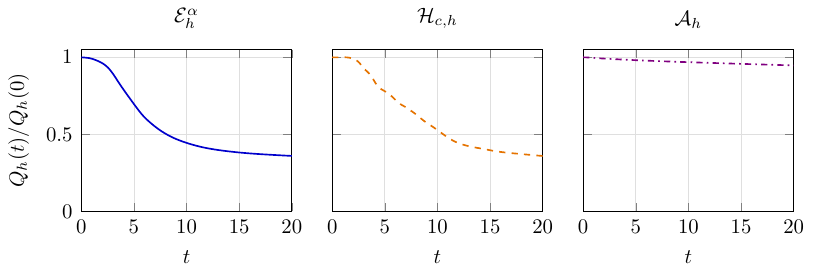}%
    \caption{Selective decay of the three invariants, each normalized by its initial value.}
    \label{fig:2d-ozt-selective-decay}
\end{figure}

We now focus on the non-ideal case and investigate the different dissipation rates of the ideal invariants. We repeat the Orszag--Tang test with $\nu = \eta = 10^{-3}$ and $S=1$ on the $128\times128$ mesh with $\alpha=2h$. \Cref{fig:2d-ozt-selective-decay} reports the time evolution of the three invariants. For these initial data and parameters, the structure-preserving scheme reproduces the separation of timescales reported in \cite{biskamp2003magnetohydrodynamic}: $\mathcal{E}^{\alpha}$ and the magnitude of $\mathcal{H}_c$ decay rapidly while $\mathcal{A}$ changes more slowly.

\begin{figure}[htpb]
    \centering
    \includegraphics{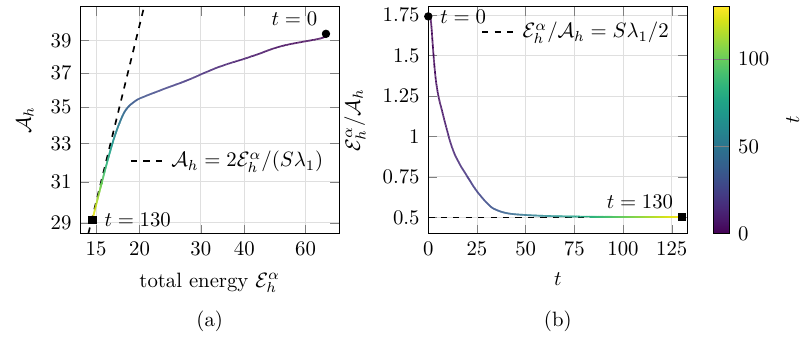}
    \caption{Orszag--Tang run up to $T=130$, colored by time. (a) Evolution of $(\mathcal{E}_h^{\alpha}(t), \mathcal{A}_h(t))$. (b) Ratio $\mathcal{E}_h^{\alpha}(t)/\mathcal{A}_h(t)$. The dashed lines are the Arnold bound $\mathcal{A}_h \le 2\mathcal{E}_h^{\alpha}/(S\lambda_1)$, i.e.\ $\mathcal{E}_h^{\alpha}/\mathcal{A}_h\ge S\lambda_1/2$, with $S=\lambda_1=1$.}
    \label{fig:2d-ozt-phase}
\end{figure}

The relaxed state is characterized by a discrete Poincar\'e-type bound. In the periodic setting, taking $e_h=1$ in \eqref{eq:sp-scheme-2d-A} shows that $\int_\Omega A_h\dx$ is conserved, so $A_h$ keeps the zero mean of the initial data. Since $\bm B_h=\curl A_h$ with $A_h\in V_h\subset H^1(\Omega)$, the Poincar\'e inequality $\|\nabla A_h\|^2\ge\lambda_1\|A_h\|^2$, with $\lambda_1=1$ the smallest nonzero eigenvalue of $-\Delta$ on $\Omega$, gives
\begin{equation}\label{eq:arnold-bound}
\mathcal{A}_h\le\frac{2\mathcal{E}_{M,h}}{S\lambda_1}\le\frac{2\mathcal{E}_h^{\alpha}}{S\lambda_1},
\qquad \mathcal{E}_{M,h}:=\frac S2\|\bm B_h\|^2.
\end{equation}
This is an analog of the three-dimensional Arnold inequality $\|\bm B\|^2\ge C|\mathcal{H}_m|$, $\mathcal{H}_m=\int_\Omega\bm A\cdot\bm B\dx$, which bounds the magnetic energy from below by the magnetic helicity \cite{ArnoldTopologicalMethodsHydrodynamics2021}; in two dimensions $\mathcal{H}_m$ vanishes and $\mathcal{A}$ takes its role. For the numerical behavior of the three-dimensional bound we refer to \cite{he2025helicity,parkerE3}. The bound is approached only when $A_h$ is concentrated in the lowest Fourier modes and the kinetic energy is small. For the original two-dimensional MHD equations on the torus, it was proposed in \cite{Matthaeus1980} and then shown in
\cite{zhan2010selective}, under a condition on the ratio of viscosity to resistivity, that
the kinetic energy decays faster than $\mathcal{A}$ and that the ratio of the energy to
$\mathcal{A}$ converges to an eigenvalue $\lambda$ of the Laplacian, the magnetic field
approaching a force-free state on which $\mathcal{A}$ decays like $e^{-2\eta\lambda t}$.
The Arnold-type bound corresponds to the lowest eigenvalue $\lambda_1$, where $\mathcal{A}$
decays like $e^{-2\eta\lambda_1 t}$. \Cref{fig:2d-ozt-phase} shows a run up to $T=130$.
In this experiment the energy decays much faster than $\mathcal{A}_h$, and we observe that the
evolution approaches the bound and then relaxes along it, which is consistent with the selective decay described in
\cite{zhan2010selective}. 



\subsection{Large-scale simulation and solver robustness test}\label{sec:ozt-large-scale}

\begin{figure}[htpb]
    \centering
    \begin{subfigure}[b]{0.45\textwidth}
        \centering
        \includegraphics[width=\linewidth]{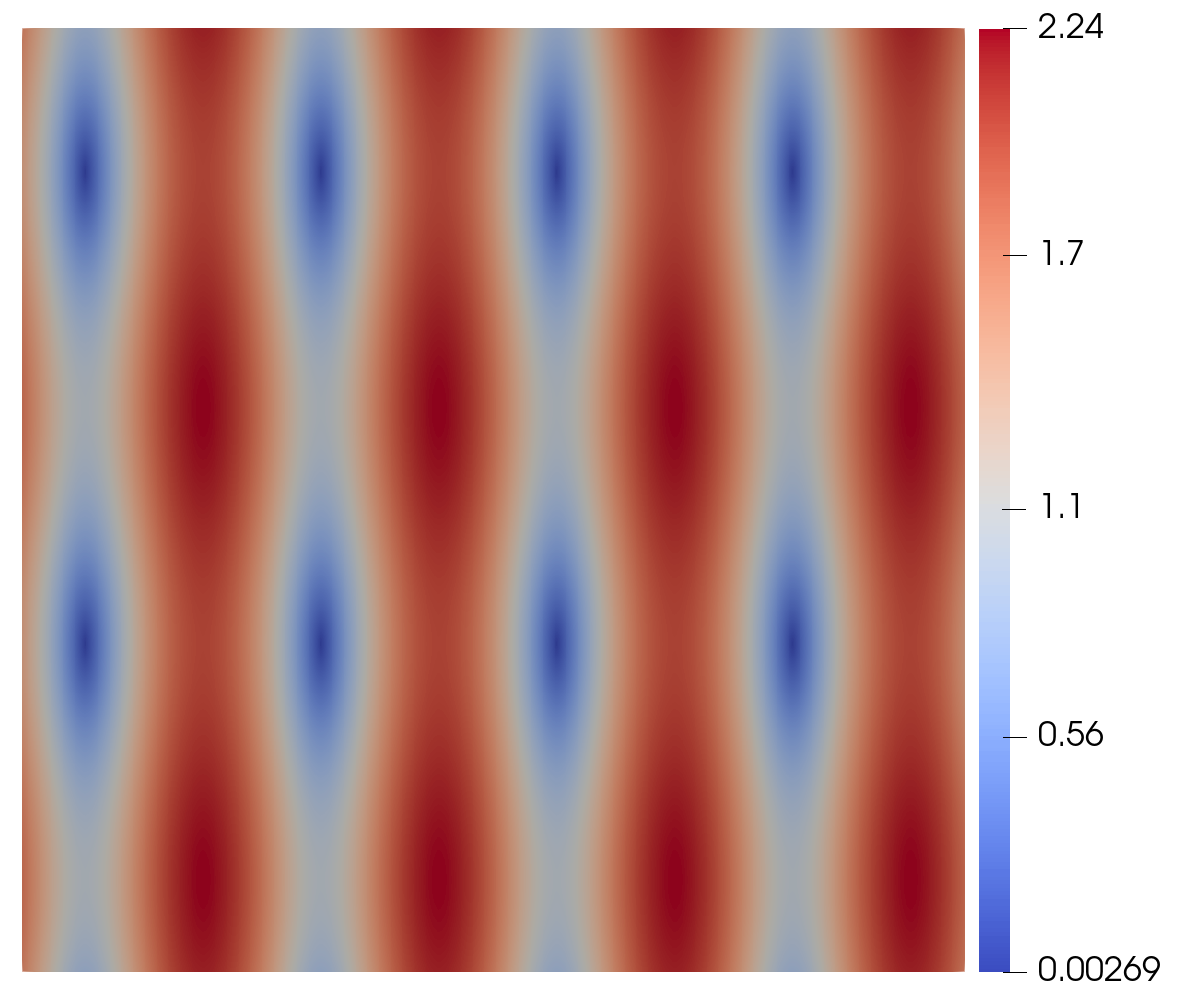}
        \caption{$t=0$}
    \end{subfigure}
    \hfill
    \begin{subfigure}[b]{0.45\textwidth}
        \centering
        \includegraphics[width=\linewidth]{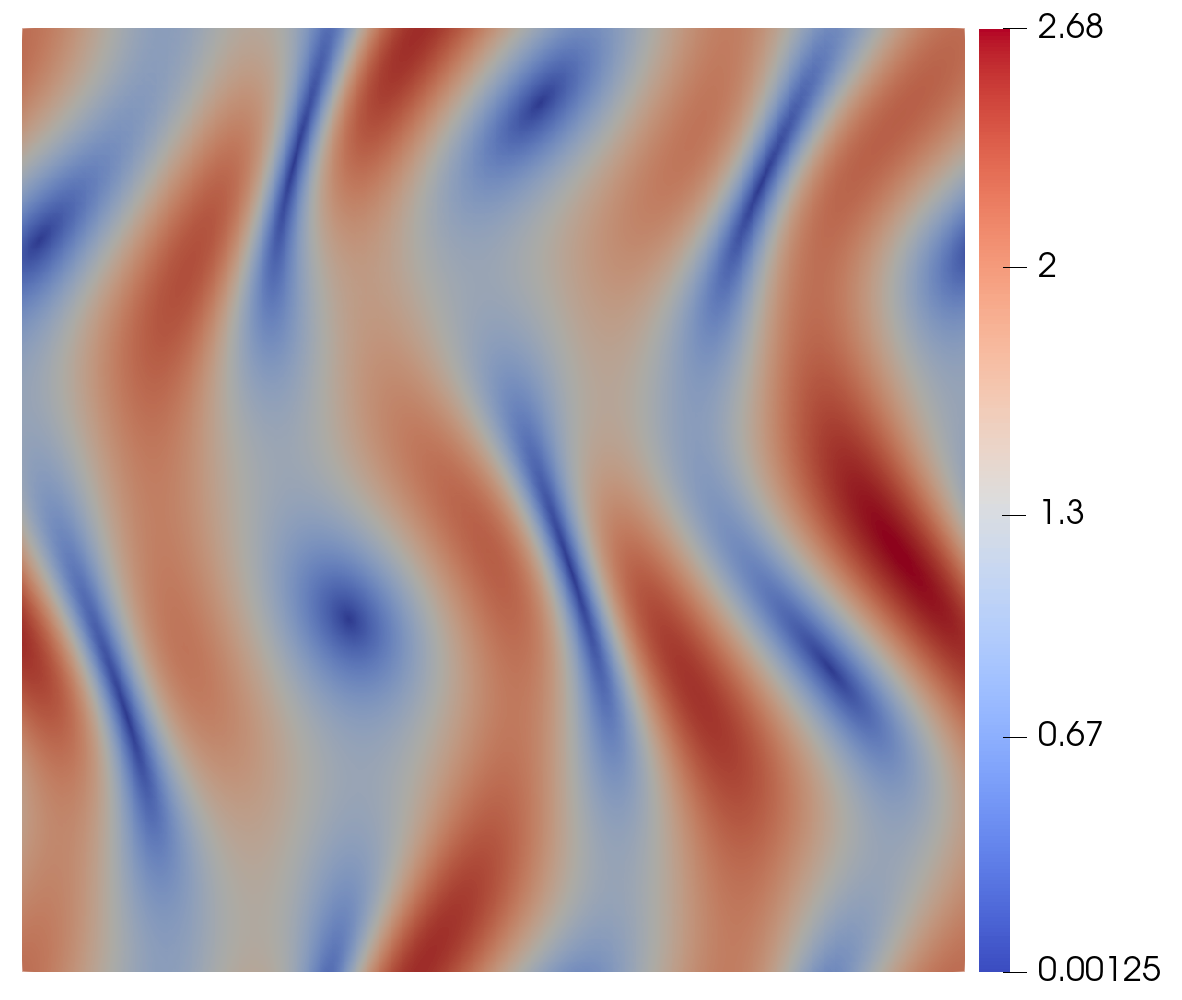}
        \caption{$t=0.5$}
    \end{subfigure}
    \begin{subfigure}[b]{0.45\textwidth}
        \centering
        \includegraphics[width=\linewidth]{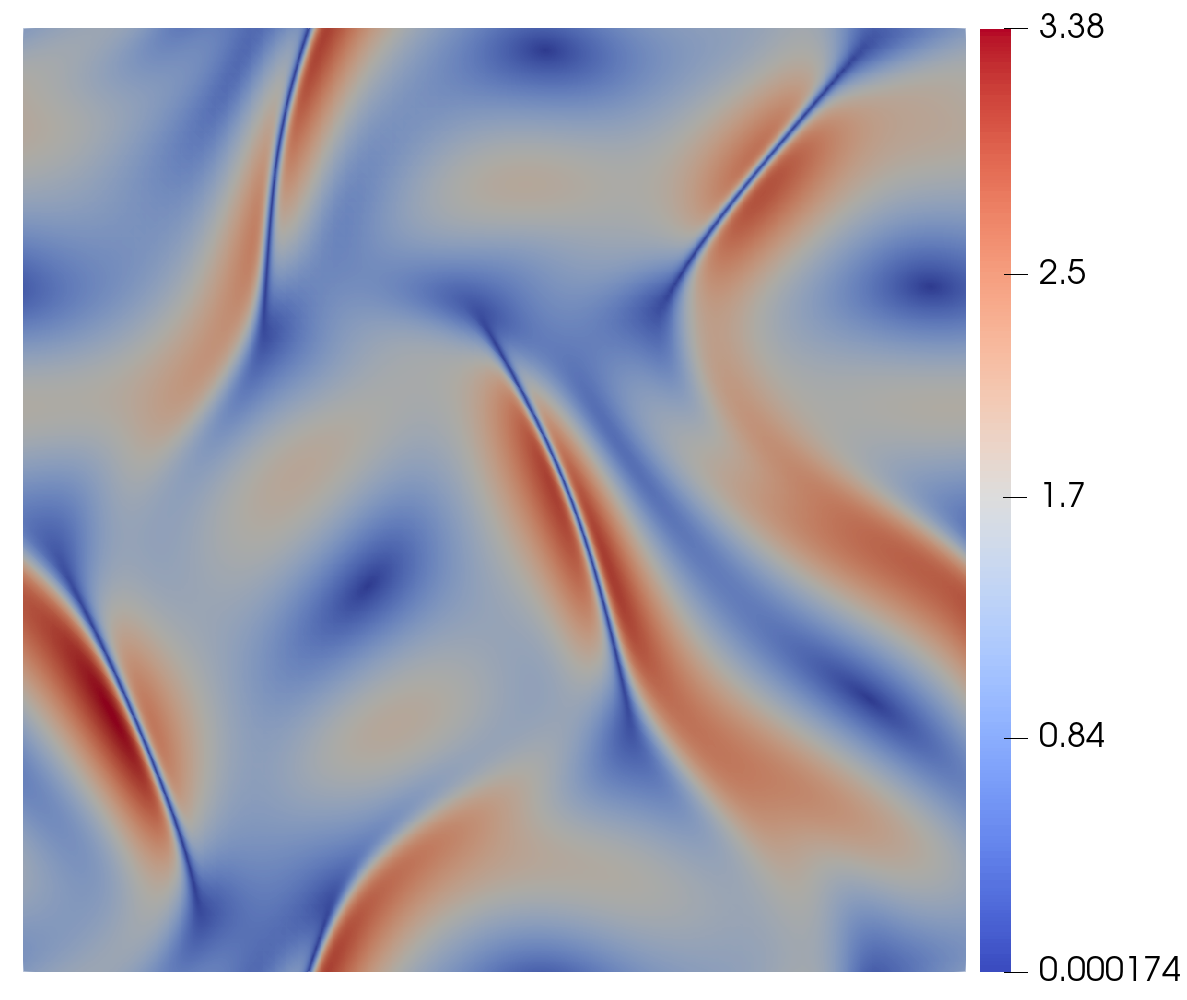}
        \caption{$t=1$}
    \end{subfigure}
    \hfill
    \begin{subfigure}[b]{0.45\textwidth}
        \centering
        \includegraphics[width=\linewidth]{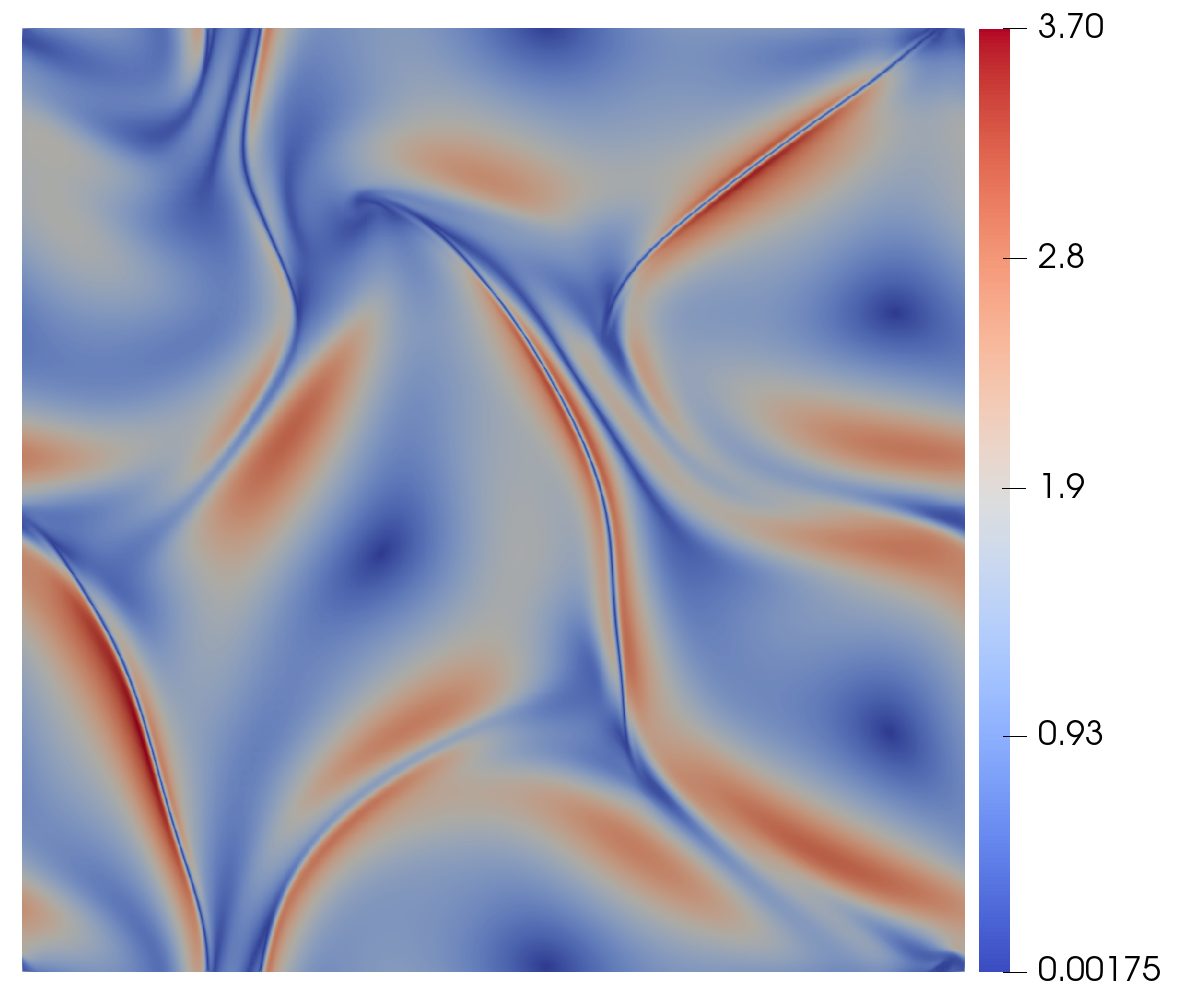}
        \caption{$t=1.5$}
    \end{subfigure}
    \begin{subfigure}[b]{0.45\textwidth}
        \centering
        \includegraphics[width=\linewidth]{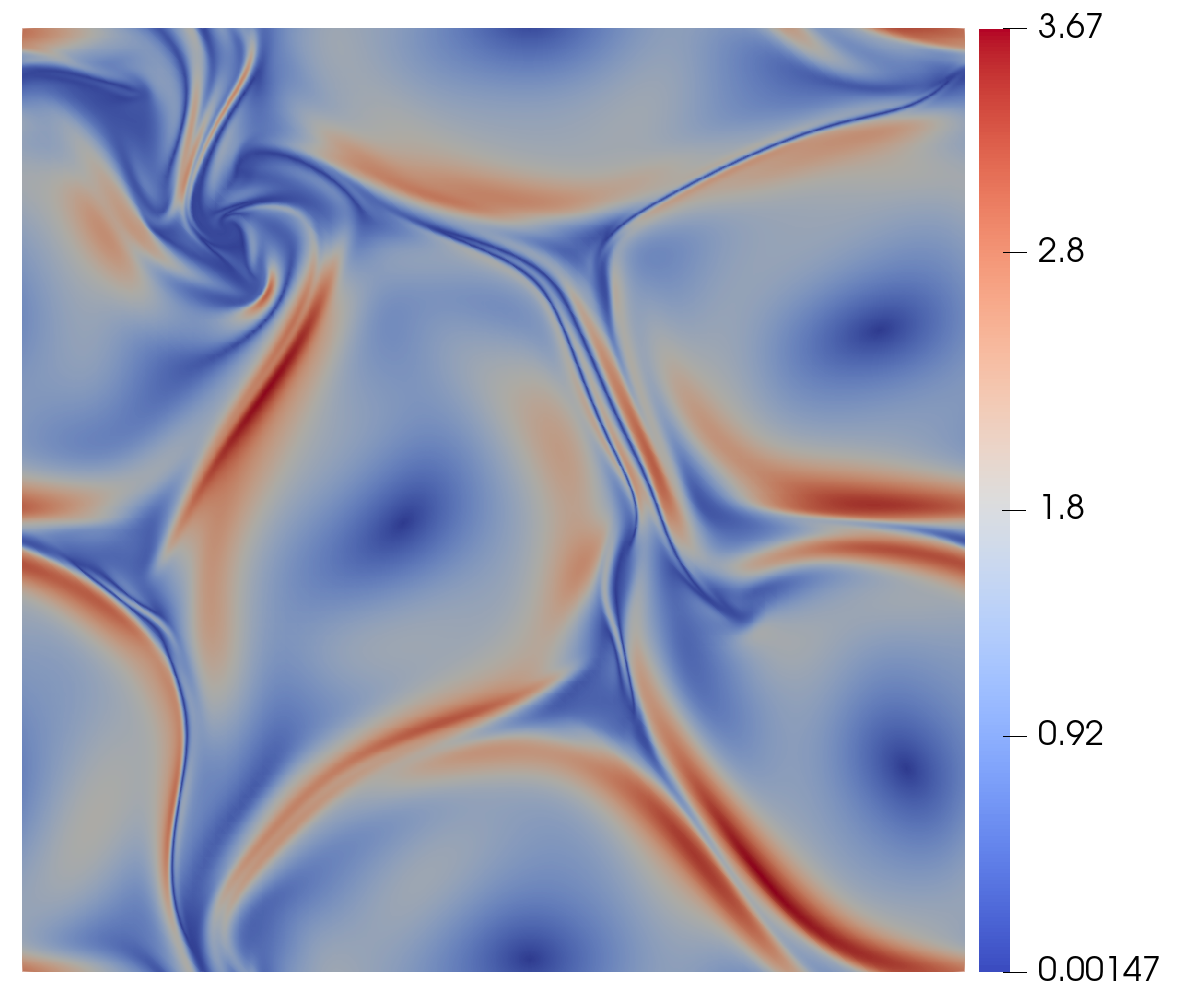}
        \caption{$t=2$}
    \end{subfigure}
    \hfill
    \begin{subfigure}[b]{0.45\textwidth}
        \centering
        \includegraphics[width=\linewidth]{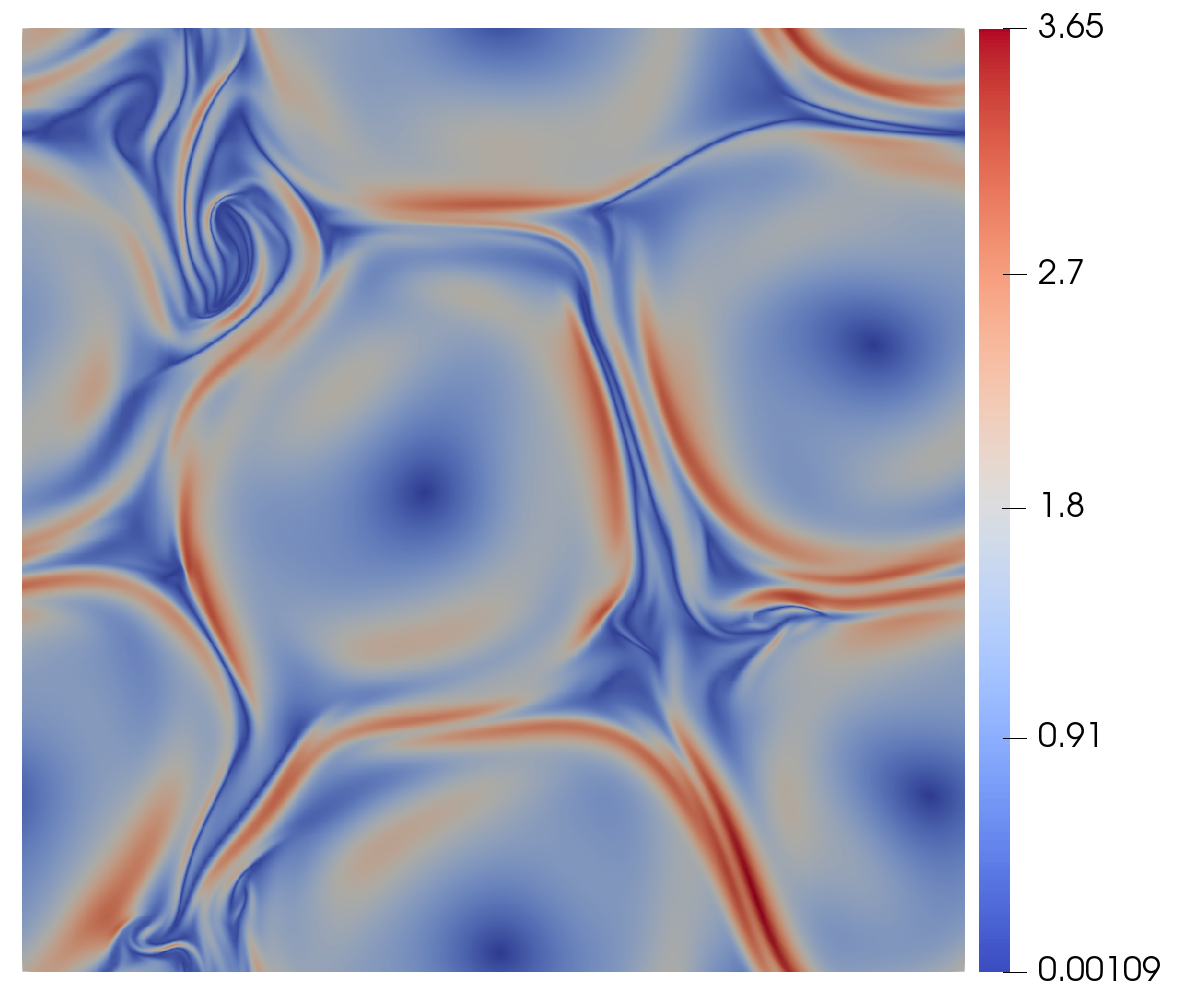}
        \caption{$t=2.5$}
    \end{subfigure}
    \caption{Orszag--Tang test on the $1024\times1024$ mesh with $\alpha=2h$ and $\Delta t=0.01$, colored by the magnetic field strength $|\bm B_h|$; the color range is rescaled in each frame.}
    \label{fig:ozt-simulation}
\end{figure}

We then run a large-scale simulation until $T=2.5$. For the multigrid hierarchy we use a coarse mesh of $16\times 16$ cells and $6$ levels of refinement, resulting in a $1024\times1024$ mesh with $2\,097\,152$ cells and $12\,582\,912$ degrees of freedom, and we take $\Delta t=0.01$. \Cref{fig:ozt-simulation} shows the evolution of the magnetic field at $t=0,0.5,\dots,2.5$. Starting from smooth, large-scale fields, the magnetic eddies steepen into thin current sheets, which reconnect and destabilize, driving a cascade to small-scale MHD turbulence.

We further test the robustness of the solver of \Cref{sec:solver-development} on the Orszag--Tang problem with $\nu=\eta=10^{-3}$, for $\alpha\in\{2h,4h\}$ and $\Delta t\in\{0.01,0.05\}$, on meshes obtained by refining a coarse mesh of $16\times16$ cells. We run until $T=2$, which covers the formation of the thinnest current sheets, and \Cref{tab:newton-performance} reports the average number of Newton iterations and the average number of outer FGMRES iterations per Newton iteration over all time steps. For the small time step $\Delta t=0.01$, the outer iteration counts do not grow under mesh refinement and are nearly independent of the filter width, while the average number of Newton steps increases mildly, from $2$ to about $2.5$. For the larger time step $\Delta t=0.05$, the iteration counts grow under mesh refinement for $N\ge256$, but they are smaller for $\alpha=4h$ than for $\alpha=2h$. This is consistent with the physics: as the current sheets thin, increasingly small-scale structures develop on finer meshes, and with a large time step each Newton solve must resolve a larger change of the solution. A larger $\alpha$ smooths these structures and therefore mitigates the growth.
\begin{table}[htbp]
\centering
\caption{Solver performance: average number of Newton steps (in parentheses) and average number of outer FGMRES iterations per Newton step, both averaged over all time steps up to $T=2$ ($200$ steps for $\Delta t=0.01$, $40$ steps for $\Delta t=0.05$) on $N\times N$ meshes.}
\label{tab:newton-performance}

\begin{tabular}{rrrrrr}
\toprule
& &
\multicolumn{2}{c}{$\Delta t = 0.01$} &
\multicolumn{2}{c}{$\Delta t = 0.05$} \\
\cmidrule(lr){3-4}
\cmidrule(lr){5-6}
$N$ & DOFs
& $\alpha=2h$ & $\alpha=4h$
& $\alpha=2h$ & $\alpha=4h$ \\
\midrule
32  & 12,288    & (2.00) 4.76 & (2.00) 4.94 & (2.60) 3.39 & (2.02) 3.49 \\
64  & 49,152    & (2.00) 4.35 & (2.00) 4.41 & (2.88) 3.16 & (2.62) 3.07 \\
128 & 196,608   & (2.15) 3.77 & (2.00) 3.79 & (3.15) 3.60 & (2.88) 3.28 \\
256 & 786,432   & (2.46) 3.56 & (2.23) 3.54 & (3.30) 4.97 & (3.17) 4.35 \\
512 & 3,145,728 & (2.53) 3.67 & (2.47) 3.59 & (3.40) 7.71 & (3.33) 6.37 \\
\bottomrule
\end{tabular}

\end{table}

\subsection{Inverse cascade for random initial data}\label{sec:random-relaxation}
The inverse cascade of $\mathcal{A}$ \cite{fyfe1976high,pouquet1978two,biskamp1994dynamics}
is the mechanism commonly invoked for the relaxation to the bound \eqref{eq:arnold-bound}
\cite{Matthaeus1980,biskamp2003magnetohydrodynamic}: we now test whether the scheme captures it starting from band-limited random data.

\textit{Initial data and parameters.} We consider the long-time relaxation of a random configuration on the periodic square $\Omega = (0, 2\pi)^{2}$ with no external forcing. We define the following wavenumber set
\begin{equation}\label{eqn:decay-shell}
\mathcal{K} := \Bigl\{\mathbf{k} = (k_x, k_y)\in \mathbb{Z}^{2}\;:\;
        k_x \ge 1\ \text{or}\ (k_x = 0,\ k_y \ge 1),\
        \bigl||\mathbf{k}| - k_{\mathrm p}\bigr| \le 2\Bigr\},
\end{equation}
with peak wavenumber $k_{\mathrm p} = 20$. The velocity stream function $\varphi$ and the
magnetic flux function $\psi$ are sampled as
\begin{equation}\label{eqn:decay-stream-flux}
\varphi(\mathbf{x}) := \sum_{\mathbf{k}\in \mathcal{K}}
       \frac{a}{|\mathbf{k}|}\,\cos \bigl(\mathbf{k}\!\cdot\!\mathbf{x} + \theta^{\varphi}_{\mathbf{k}}\bigr),
\qquad
\psi(\mathbf{x}) := \sum_{\mathbf{k}\in \mathcal{K}}
       \frac{a}{|\mathbf{k}|}\,\cos\bigl(\mathbf{k}\!\cdot\!\mathbf{x} + \theta^{\psi}_{\mathbf{k}}\bigr),
\end{equation}
with amplitude $a = 0.5$ and i.i.d.\ random phases
$\theta^{\varphi}_{\mathbf{k}},\,\theta^{\psi}_{\mathbf{k}} \sim
\mathcal{U}(0, 2\pi)$, drawn once with NumPy's
\texttt{default\_rng(42)} so the realization is reproducible: the modes are visited with $k_x$ and then $k_y$ in
increasing order, and for each mode $\theta^{\psi}_{\mathbf{k}}$ is drawn before $\theta^{\varphi}_{\mathbf{k}}$.

The functions $\varphi$ and $\psi$ are interpolated into $V_h$, and the initial conditions are
\begin{equation}\label{eqn:decay-init}
\bm u_h^0 = \nabla^{\perp}\!\varphi_h, \qquad A_h^0 = -\psi_h, \qquad \bm B_h^0 = \curl A_h^0 = \nabla^{\perp}\!\psi_h,
\end{equation}
so that $\bm u_h^0,\bm B_h^0\in\bm F_h$ are exactly
divergence-free; the filtered velocity $\ov{\bm u}_h^0$ is obtained from $\bm u_h^0$ by the discrete filter.

We take $\alpha = 2h$, a time step $\Delta t = 0.05$, and integrate to $T = 75$ on a $128 \times 128$ mesh
with lowest-order $k=1$ elements in space, $\nu = \eta = 10^{-3}$ and $S=1$; the nonlinear systems are solved by Newton's
method with a sparse direct LU factorization (MUMPS).

\begin{figure}[htpb]
    \centering
    \includegraphics[width=1.0\linewidth]{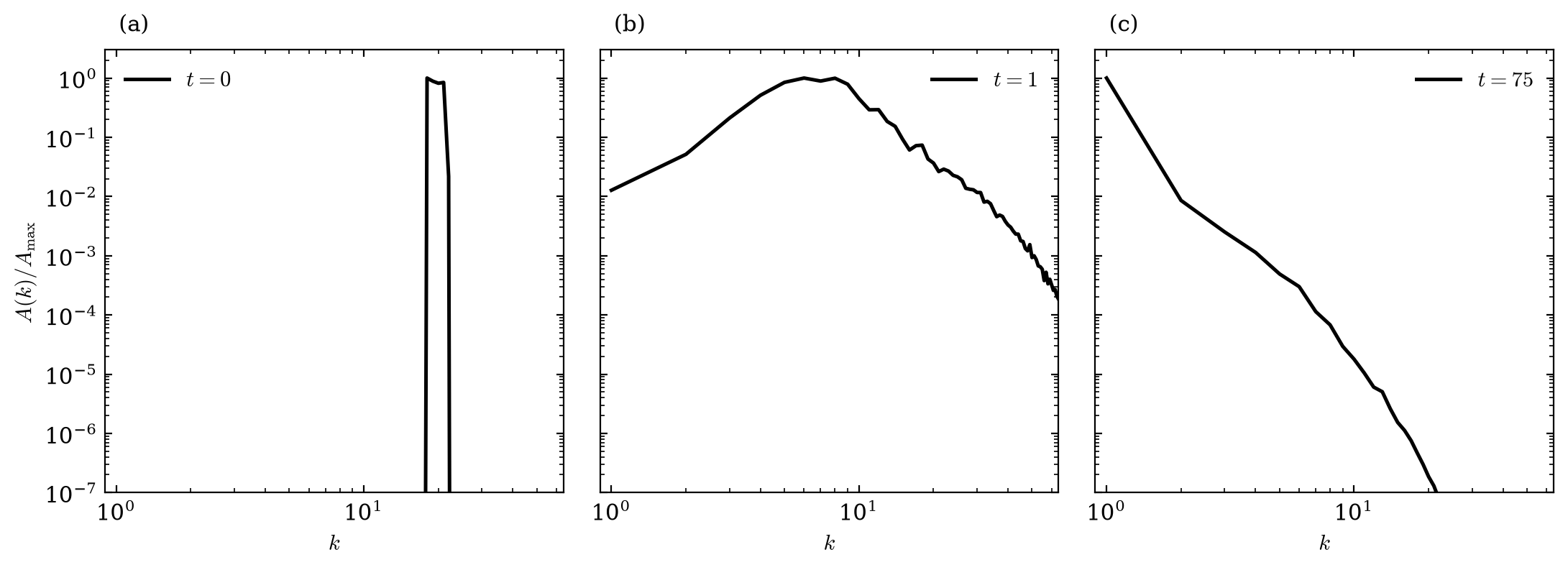}
    \caption{Mean square potential spectrum $A(k)/A_{\max}$ at $t = 0, 1, 75$, normalized at each time by its maximum (see \ref{app:spectra}). The band-limited initial state at $k\approx 20$ develops a broad spectrum whose peak migrates to $k\approx 6$ at $t=1$ and ultimately condenses 
at the box mode $k=1$ at $t=75$, demonstrating the inverse cascade of
$\mathcal{A}$.}
    \label{fig:A2-inverse-cascade}
\end{figure}

\begin{figure}[htbp]
    \centering
    \includegraphics[width=1.0\linewidth]{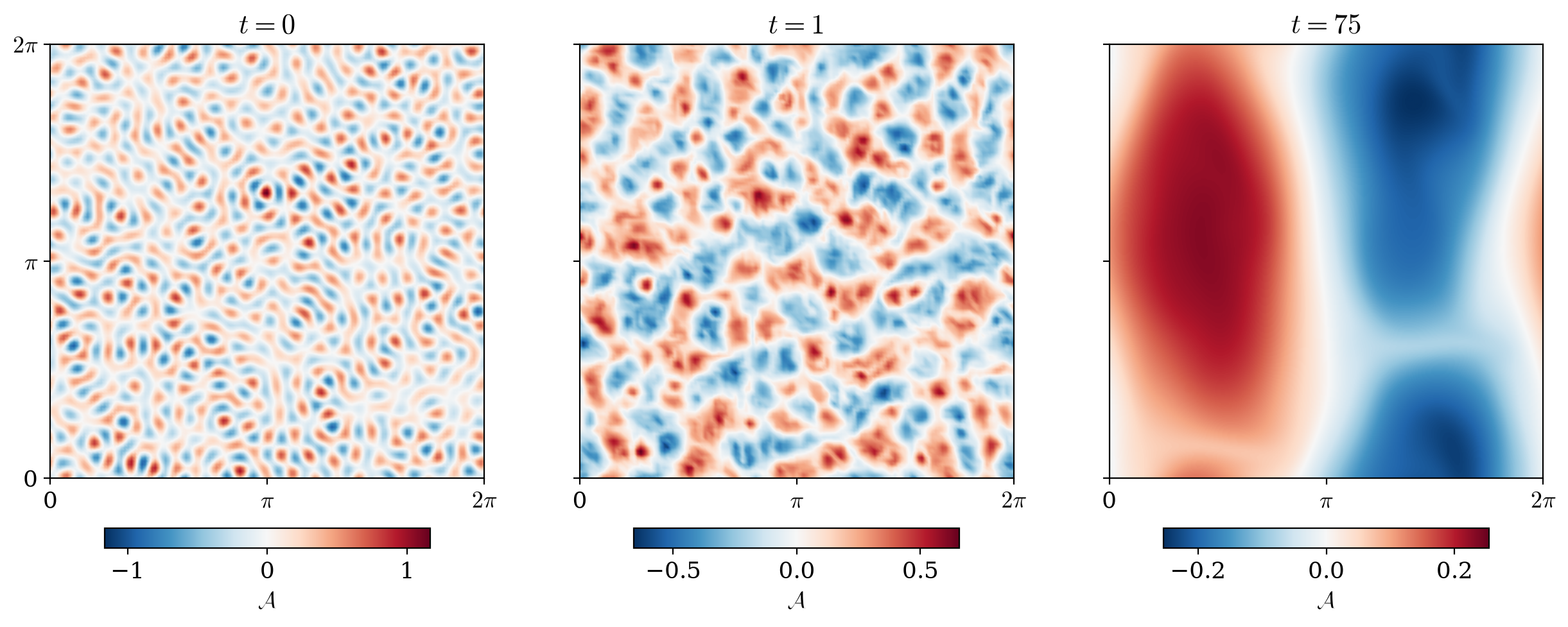}
    \caption{Evolution of the magnetic potential $A_h$.}
    \label{fig:2d-decay-A}
\end{figure}

\Cref{fig:A2-inverse-cascade} shows the inverse cascade of the mean square potential: the peak of its spectrum
moves from the initial band $k\approx20$ to the box mode $k=1$. \Cref{fig:2d-decay-A} shows the corresponding
formation of large-scale magnetic structures from the initial small-scale field.

\begin{figure}[htbp]
    \centering
    \includegraphics{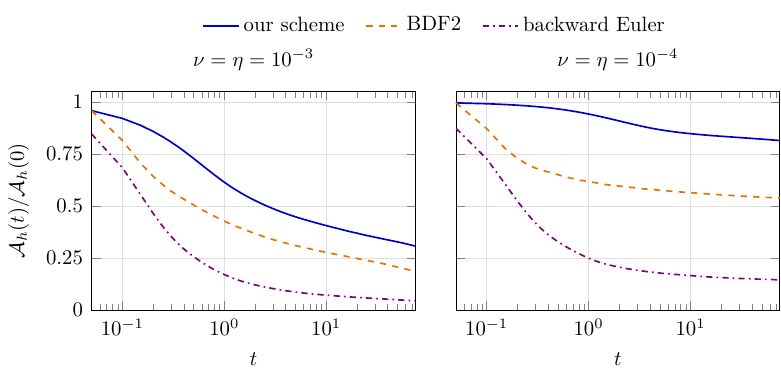}
    \caption{Mean square potential, relative to its initial value, for the random initial data on the $128\times128$ mesh with $\alpha=2h$ and $\Delta t=0.05$: the scheme of \Cref{alg:sp-discretisation-mhd-alpha}, and BDF2 and backward Euler with the same spatial discretization. Left: $\nu=\eta=10^{-3}$, the computation of \Cref{fig:A2-inverse-cascade,fig:2d-decay-A}. Right: $\nu=\eta=10^{-4}$.}
    \label{fig:2d-decay-schemes}
\end{figure}

\textit{Comparison with non-structure-preserving schemes.} To see what the discrete laws of \Cref{thm:2d-sp} contribute to this result, we repeat
the computation with the same spatial discretization and time step, replacing the
implicit midpoint rule by the BDF2 formula (with a midpoint first step) and by the
backward Euler method, neither of which satisfies the three discrete laws.
\Cref{fig:2d-decay-schemes} shows the mean square potential for $\nu=\eta=10^{-3}$ and
$\nu=\eta=10^{-4}$. Both alternative schemes lose substantially more $\mathcal{A}_h$ than
ours, backward Euler most of all, and the gap does not shrink when the dissipation is
reduced tenfold. This additional loss is therefore numerical. It occurs mostly for
$t\lesssim1$, while the small scales are active, and it changes the state toward which
the solution relaxes.

\begin{figure}[htbp]
    \centering
    \includegraphics[width=0.8\linewidth]{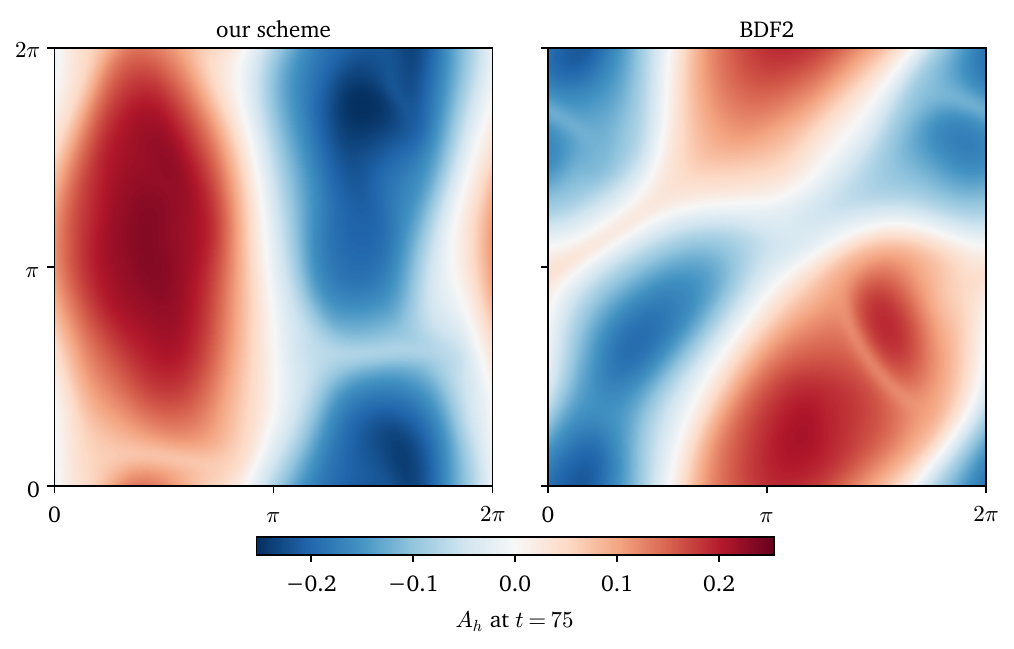}
    \caption{Magnetic potential $A_h$ (mean removed) at $t=75$ for $\nu=\eta=10^{-3}$: the scheme of \Cref{alg:sp-discretisation-mhd-alpha} (left) and BDF2 (right), with one color scale.}
    \label{fig:2d-decay-bdf2}
\end{figure}

\textit{Relaxed state comparison.} \Cref{fig:2d-decay-bdf2} compares the magnetic potential computed by our scheme and by
BDF2 at $t=75$ for $\nu=\eta=10^{-3}$. The solution of our scheme is concentrated in the first eigenspace of
$-\Delta$, the Fourier modes with $|\bm k|^2=\lambda_1=1$, and is therefore close to the bound
\eqref{eq:arnold-bound}. With BDF2 a large part of $\mathcal{A}_h$ remains in
the modes with $|\bm k|^2=2$, which appear as the diagonal structures in
\Cref{fig:2d-decay-bdf2}, and the solution is still far from the bound. 

\section{Energy spectra and the effect of filtering}\label{sec:spectra}
In two-dimensional MHD, energy cascades forward to small scales, and following it down to the dissipation scale is what makes direct simulation expensive.
In the $\alpha$-model, the term $\alpha^2\|\nabla\ov{\bm u}\|^2$ in the energy penalizes scales below $\alpha$, so the kinetic energy spectrum is expected to fall off for $k\alpha\gg1$ while the large scales are left unchanged \cite{foias2001navier,chen1999direct}. To check this, we compare kinetic energy spectra of the scheme with a Fourier pseudo-spectral direct numerical simulation of the original MHD equations computed with Dedalus \cite{burns2020dedalus}.

We use the Orszag--Tang test of \Cref{sec:selective-decay} with $\nu=\eta=10^{-4}$, $S=1$, $\Delta t=0.01$ and $\alpha\in\{0,2h,4h\}$. For the multigrid method we use a coarse mesh of $8\times 8$ cells and $5$ levels of refinement, resulting in a $256\times256$ mesh with lowest-order $k=1$ elements and $786\,432$ degrees of freedom. At this Reynolds number a pseudo-spectral simulation on a $256\times256$ grid would not resolve the dissipation scale, so the reference is computed on a $1024\times1024$ grid. The shell-summed spectrum $E_K(k)$ is the kinetic energy of all Fourier modes with
$|\bm k|\approx k$, so that $(2\pi)^2\sum_k E_K(k)=\mathcal{E}_K$. For the scheme we plot the shell-summed spectrum of the kinetic part of the energy, $\mathcal{E}_{K,h}^{\alpha}:=\frac12(\bm u_h,\ov{\bm u}_h)$, computed as described in \ref{app:spectra}, and for the reference that of $\mathcal{E}_K:=\frac12\|\bm u\|^2$. They coincide when $\alpha=0$.

\begin{figure}[htbp]
    \centering
    \includegraphics{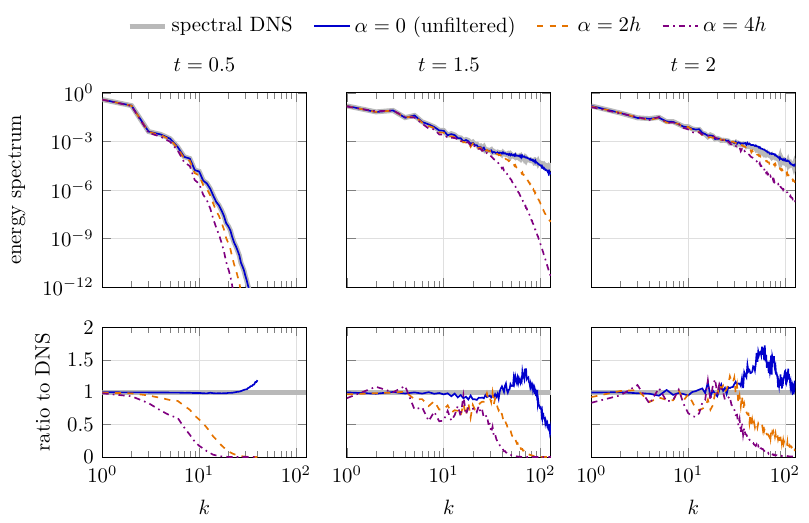}%
    \caption{Kinetic energy spectra of the Orszag--Tang test with $\nu=\eta=10^{-4}$ on the $256\times256$ mesh at
    $t=0.5$, $1.5$ and $2$ (top), and their ratio to the spectral DNS of the original MHD equations on a
    $1024\times1024$ grid (bottom). At $t=0.5$ the ratio is shown for $k\le40$.}
    \label{fig:2d-ozt-spectrum-nu1e-4}
\end{figure}

\begin{figure}[htbp]
    \centering
    \includegraphics[width=\linewidth]{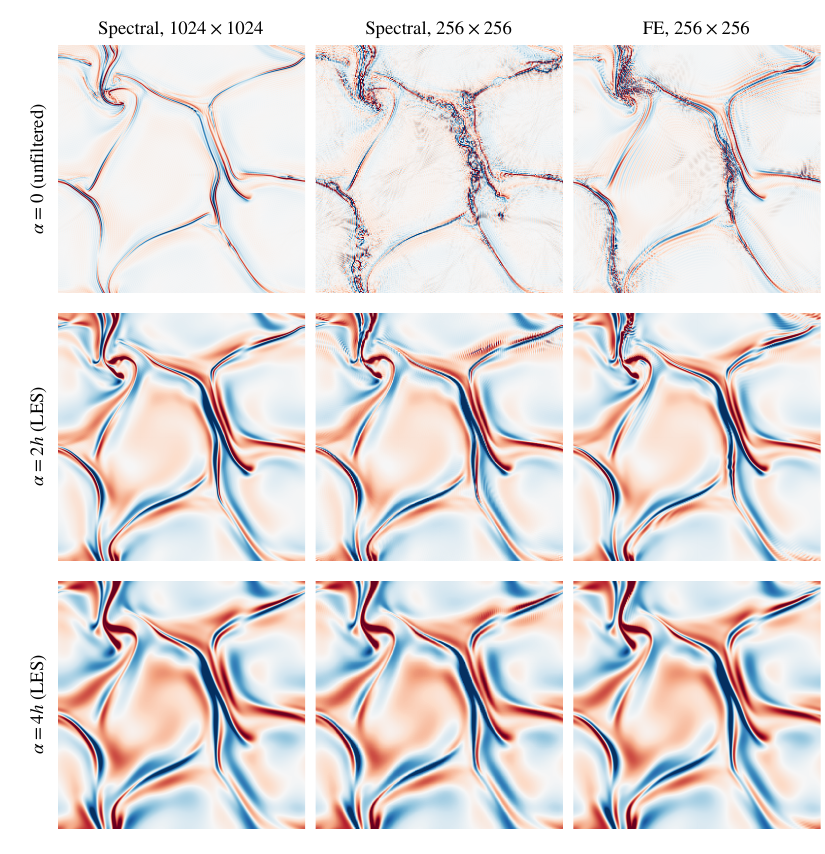}
    \caption{Vorticity of the transport velocity at $t=2$ for $\nu=\eta=10^{-4}$.
    Columns: pseudo-spectral solutions on $1024\times1024$ and $256\times256$ grids, and the
    scheme on the $256\times256$ mesh. Rows: $\alpha=0$, $2h$ and $4h$, where $h$ is the
    mesh size of the $256\times256$ mesh. The $1024\times1024$ fields are restricted to
    $|\bm k|\le128$. Each row shares one color scale.}
    \label{fig:2d-ozt-pileup-fields}
\end{figure}

\Cref{fig:2d-ozt-spectrum-nu1e-4} shows the spectra and their ratio to the reference, and
\Cref{fig:2d-ozt-pileup-fields} the corresponding vorticity fields.

\textit{Spectral space.} At large scales, all runs agree with the reference. At $t=0.5$
the run with $\alpha=0$ follows the reference over the whole range. At $t=1.5$ and $2$ the forward cascade has reached the mesh scale. The scheme has no numerical dissipation, so the energy that arrives at the smallest
resolved scales is not removed there. For $\alpha=0$ it accumulates near the grid scale, and
the spectrum rises above the reference for $k\gtrsim30$. For $\alpha>0$ the spectrum falls
off at high wavenumbers, with the cut-off moving to smaller $k$ as $\alpha$ increases, and
with $\alpha=2h$ and $4h$ the accumulation is suppressed. Thus $\alpha$ leaves the largest
scales ($k\alpha \ll 1$) unchanged and, as expected, also reduces the spectrum at
intermediate wavenumbers $k\alpha\sim1$, giving a more compact energy spectrum, in agreement
with the regularization effect of the Navier--Stokes-$\alpha$ model in
\cite{chen1999direct,foias2001navier}.

\textit{Physical space.} The accumulated energy appears as grid-scale oscillations along
the vortex sheets for $\alpha=0$; with $\alpha=2h$ and $4h$ the oscillations are absent.
\Cref{fig:2d-ozt-pileup-fields} also shows pseudo-spectral solutions of the same equations:
on a $256\times256$ grid the unfiltered pseudo-spectral solution is polluted in the same way,
while for $\alpha\in\{2h,4h\}$ both solutions on $256\times256$ cells are close to the
pseudo-spectral solution of the MHD-$\alpha$ model on a $1024\times1024$ grid. In these
experiments the MHD-$\alpha$ model therefore suppresses small-scale velocity activity below
$\alpha$, which suggests a possible reduction of the resolution requirements of a direct
simulation \cite{foias2001navier,chen1999direct}. We do not
quantify this saving here, since the mesh is the same for all values of $\alpha$ in our
experiments.

This comparison also shows why structure preservation matters when assessing a turbulence
model. Since our scheme preserves the ideal invariants and their dissipation laws
(\Cref{thm:2d-sp}), it adds no numerical dissipation, so the departure of the spectra
from the run with $\alpha=0$ can be attributed to the filter rather than to numerical
dissipation.




\section{Island coalescence}\label{sec:island}
We finally consider a two-dimensional island coalescence problem to demonstrate the effectiveness of our scheme for a physically relevant problem. We adopt the experimental setting from \cite{maoIncompressibilityDivB0Preserving2025}, in which two magnetic islands embedded in a Harris current sheet are attracted to each other by the Lorentz force and reconnect through a thin current sheet forming between them. The computational domain is the rectangle $\Omega=(0,2)\times(-\tfrac12,\tfrac12)$ and the initial state is the Fadeev equilibrium perturbed by a coalescence seed,
\begin{align}\label{eq:island-ic}
A_0(x,y)&=\delta\ln\!\Big(\cosh\frac{y}{\delta}+\varepsilon\cos\frac{x}{\delta}\Big)+\kappa\cos(\pi x)\cos\!\Big(\frac{\pi y}{2}\Big),\qquad\\
\bm B_0&=\curl A_0,\qquad \\
\bm u_0&=\bm 0,
\end{align}
with $\delta=1/(2\pi)$, $\varepsilon=0.2$ and $\kappa=0.01$, so that the unperturbed equilibrium has its O-points at $x=1/2$ and $x=3/2$ and the X-point at $x=1$. We take $S=1$ and $\mathrm{Re}=\mathrm{Re}_m=10^{5}$. Following \cite{maoIncompressibilityDivB0Preserving2025}, the walls are slip walls with the tangential magnetic field held at its initial value,
\begin{equation}\label{eq:island-bc}
\bm u\cdot\bm n=\ov{\bm u}\cdot\bm n=0,\qquad w=\ov w=0,\qquad \bm B\cdot\bm\tau=\bm B_0\cdot\bm\tau\qquad\text{on }\partial\Omega,
\end{equation}
where $\bm\tau=(-n_2,n_1)$, so that the last condition is $\bm B\times\bm n=\bm B_0\times\bm n$. The spaces of \Cref{ssec:spaces} are taken on a triangulation of $\Omega$, and we write
\begin{equation}
\bm F_{h,0}:=\bm F_h\cap\bm H_0(\div,\Omega),\qquad V_{h,0}:=V_h\cap H^1_0(\Omega).
\end{equation}
The velocity conditions in \eqref{eq:island-bc} are built into the spaces of $\bm u_h,\ov{\bm u}_h\in\bm F_{h,0}$
and $w_h,\ov w_h\in V_{h,0}$. The magnetic variables $A_h,j_h\in V_h$ carry no essential
condition, and $\bm B_h^n:=\curl A_h^n$ as in \Cref{sec:turbulence-scheme}; the tangential datum $\bm B_0\cdot\bm\tau$ enters the current equation as a natural boundary term.

The full scheme reads as follows:
given $\bm u_h^0=\bm 0$ and $A_h^0=I_hA_0$, set $\ov{\bm u}_h^0:=\mathcal F_\alpha\bm u_h^0=\bm 0$.
For $n=0,1,\ldots,N-1$, find
\begin{equation}
\begin{gathered}
(\bm u_h^{n+1},\,\ov{\bm u}_h^{n+1},\,P_h^{n+1/2},\,w_h^{n+1/2},\,\ov w_h^{n+1},\,A_h^{n+1},\,j_h^{n+1/2})\\
\in\bm F_{h,0}\times\bm F_{h,0}\times Q_h\times V_{h,0}\times V_{h,0}\times V_h\times V_h
\end{gathered}
\end{equation}
such that for all test functions $(\bm v_h,\bm\chi_h,q_h,z_h,\xi_h,e_h,k_h)$ in the same spaces:
\begin{subequations}\label{eq:island-scheme}
\begin{align}
&\frac{1}{\Delta t}(\bm u_h^{n+1}-\bm u_h^n,\bm v_h)-(\ov{\bm u}_h^{n+1/2}\times w_h^{n+1/2},\bm v_h)
 +\nu(\curl w_h^{n+1/2},\bm v_h)\notag\\
&\qquad-(P_h^{n+1/2},\div\bm v_h)-S(j_h^{n+1/2}\times\bm B_h^{n+1/2},\bm v_h)=0,\\
&(\ov{\bm u}_h^{n+1},\bm\chi_h)+\alpha^2(\curl\ov w_h^{n+1},\bm\chi_h)=(\bm u_h^{n+1},\bm\chi_h),\\
&(\div\bm u_h^{n+1},q_h)=0,\\
&(w_h^{n+1/2},z_h)=(\bm u_h^{n+1/2},\curl z_h),\\
&(\ov w_h^{n+1},\xi_h)=(\ov{\bm u}_h^{n+1},\curl\xi_h),\\
&\frac{1}{\Delta t}(A_h^{n+1}-A_h^n,e_h)+\eta(j_h^{n+1/2},e_h)-(\ov{\bm u}_h^{n+1/2}\times\bm B_h^{n+1/2},e_h)=0,\\
&(j_h^{n+1/2},k_h)=(\bm B_h^{n+1/2},\curl k_h)+\langle\bm B_0\cdot\bm\tau,k_h\rangle_{\partial\Omega}.
 \label{eq:island-current}
\end{align}
\end{subequations}

Compared with \Cref{alg:sp-discretisation-mhd-alpha}, only the spaces and the boundary term in the current
equation change.
The scheme still preserves $\div\bm u_h=\div\ov{\bm u}_h=\div\bm B_h=0$ pointwise. The boundary conditions allow nonzero normal magnetic flux, and the boundary terms in the cross-helicity balance need not vanish, so the cross helicity is not an invariant of this configuration and we do not report it here.

\begin{figure}[htbp]
    \centering
    \begin{subfigure}[b]{0.48\textwidth}
        \centering
        \includegraphics[width=\linewidth]{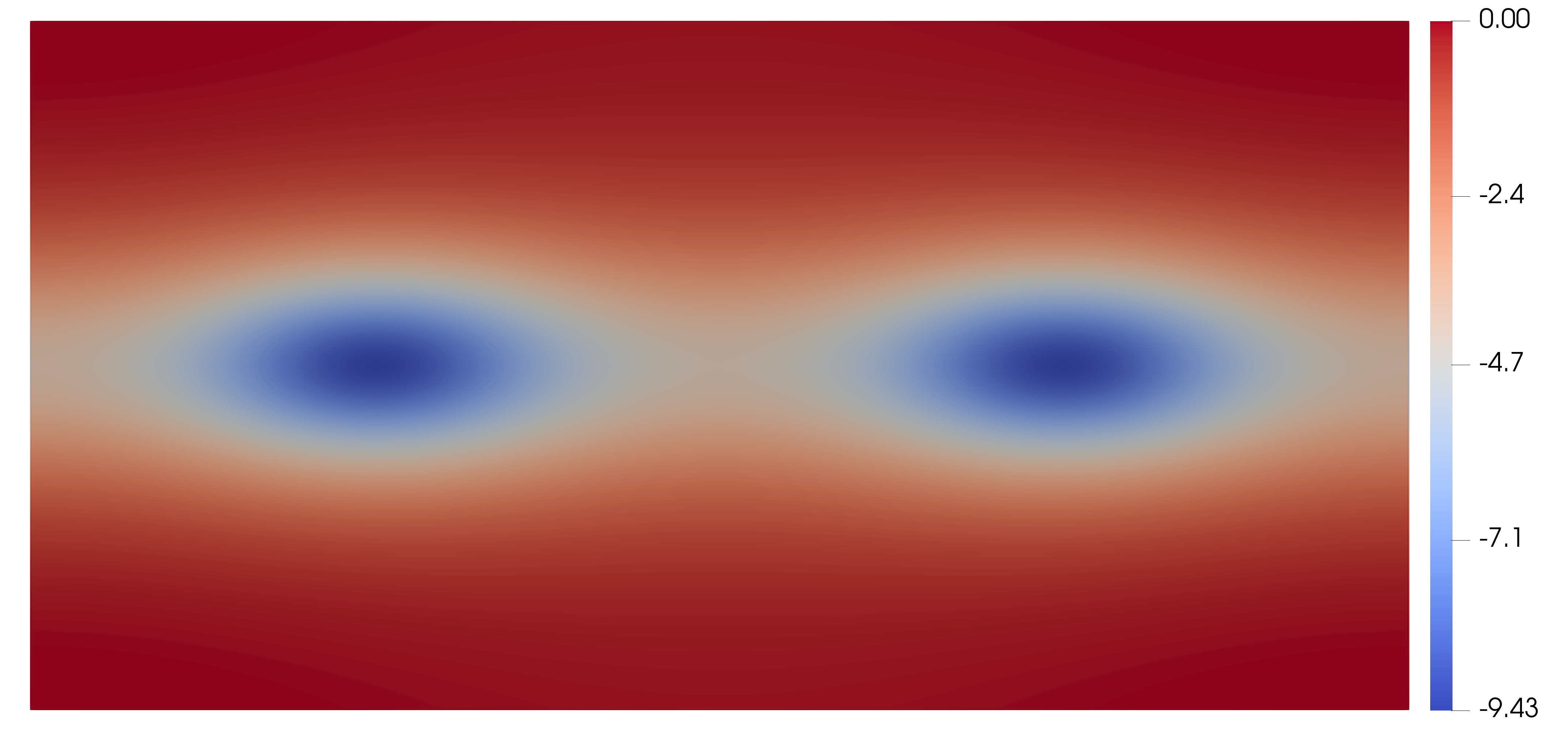}
        \caption{$t=0$}
    \end{subfigure}\hfill
    \begin{subfigure}[b]{0.48\textwidth}
        \centering
        \includegraphics[width=\linewidth]{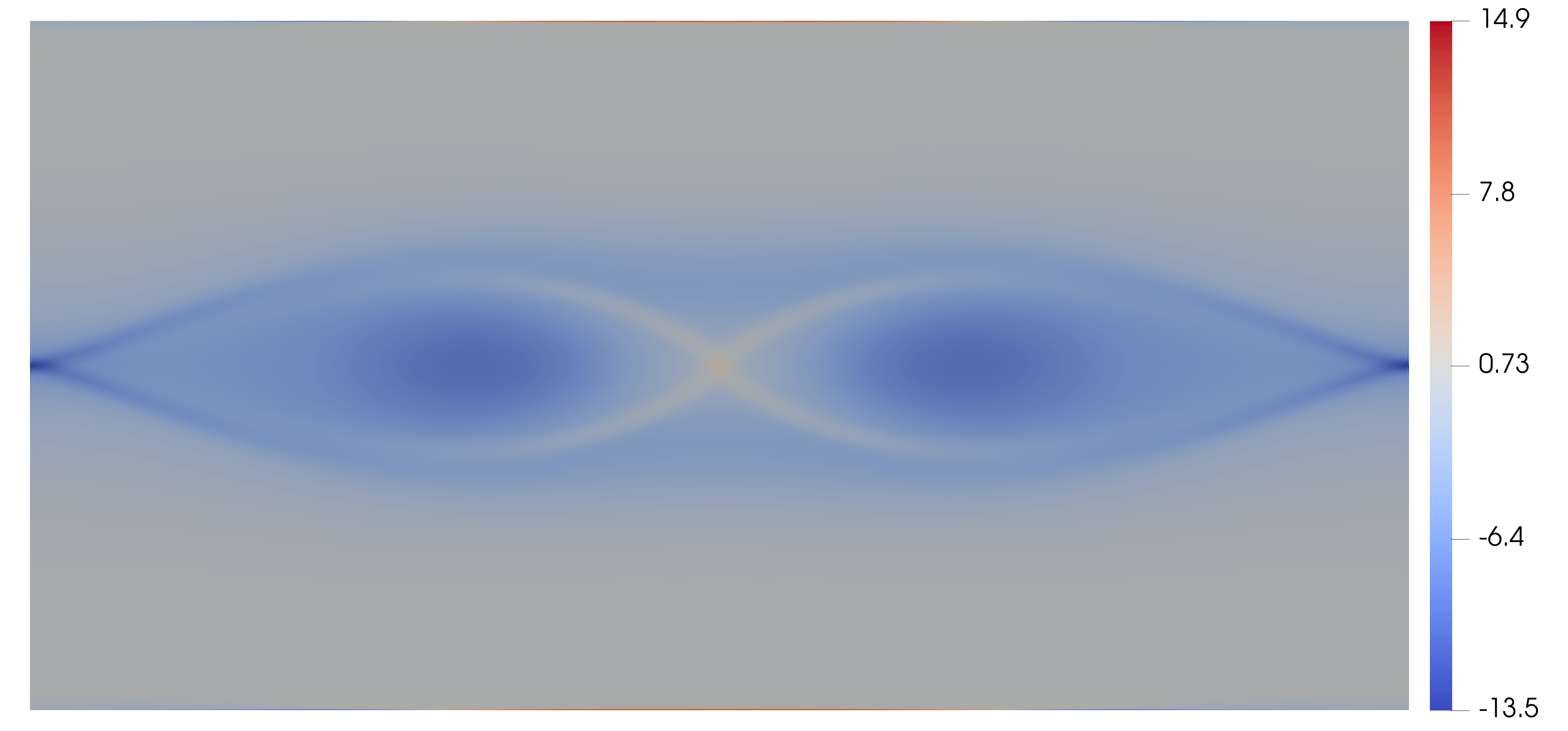}
        \caption{$t=1.5$}
    \end{subfigure}\\[1ex]
    \begin{subfigure}[b]{0.48\textwidth}
        \centering
        \includegraphics[width=\linewidth]{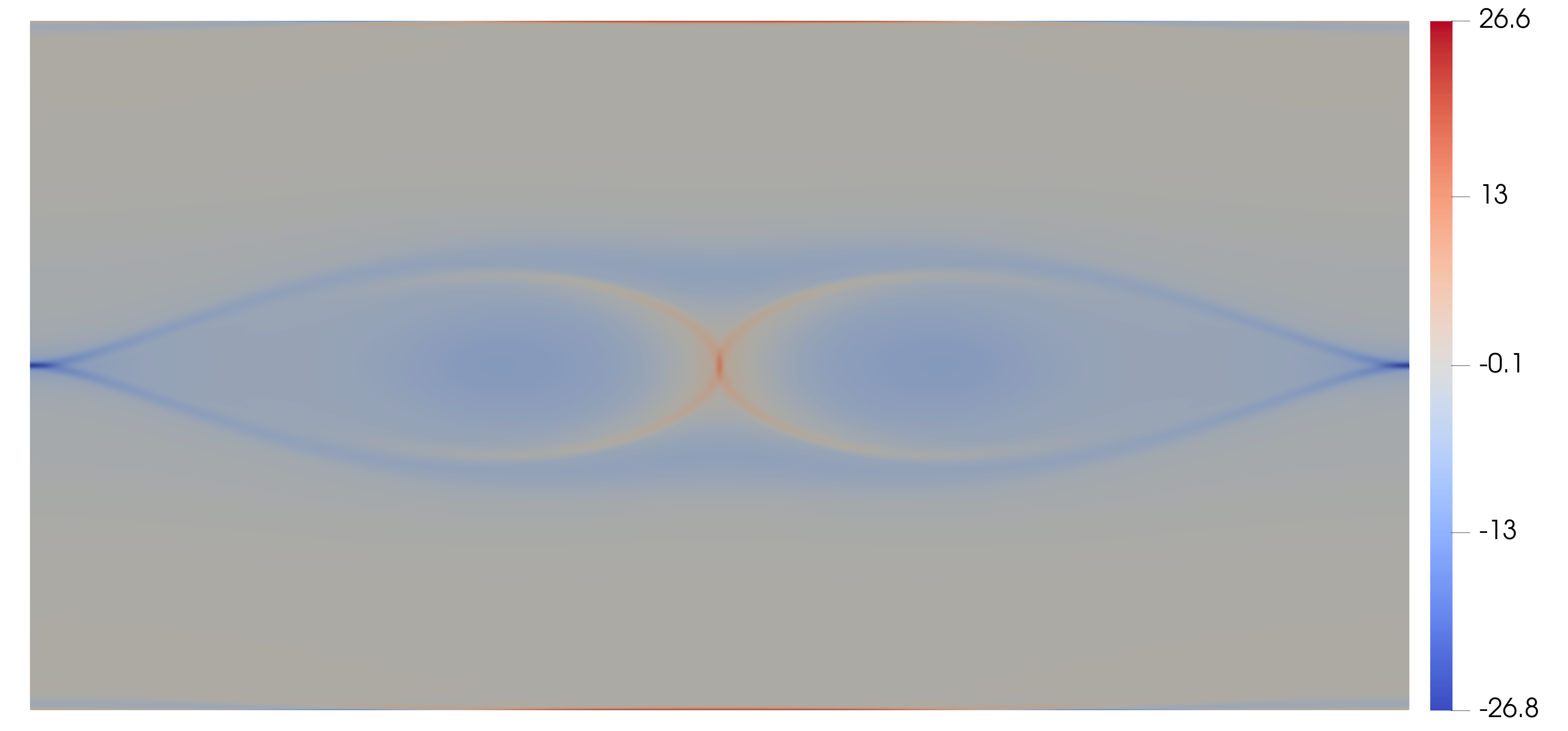}
        \caption{$t=2$}
    \end{subfigure}\hfill
    \begin{subfigure}[b]{0.48\textwidth}
        \centering
        \includegraphics[width=\linewidth]{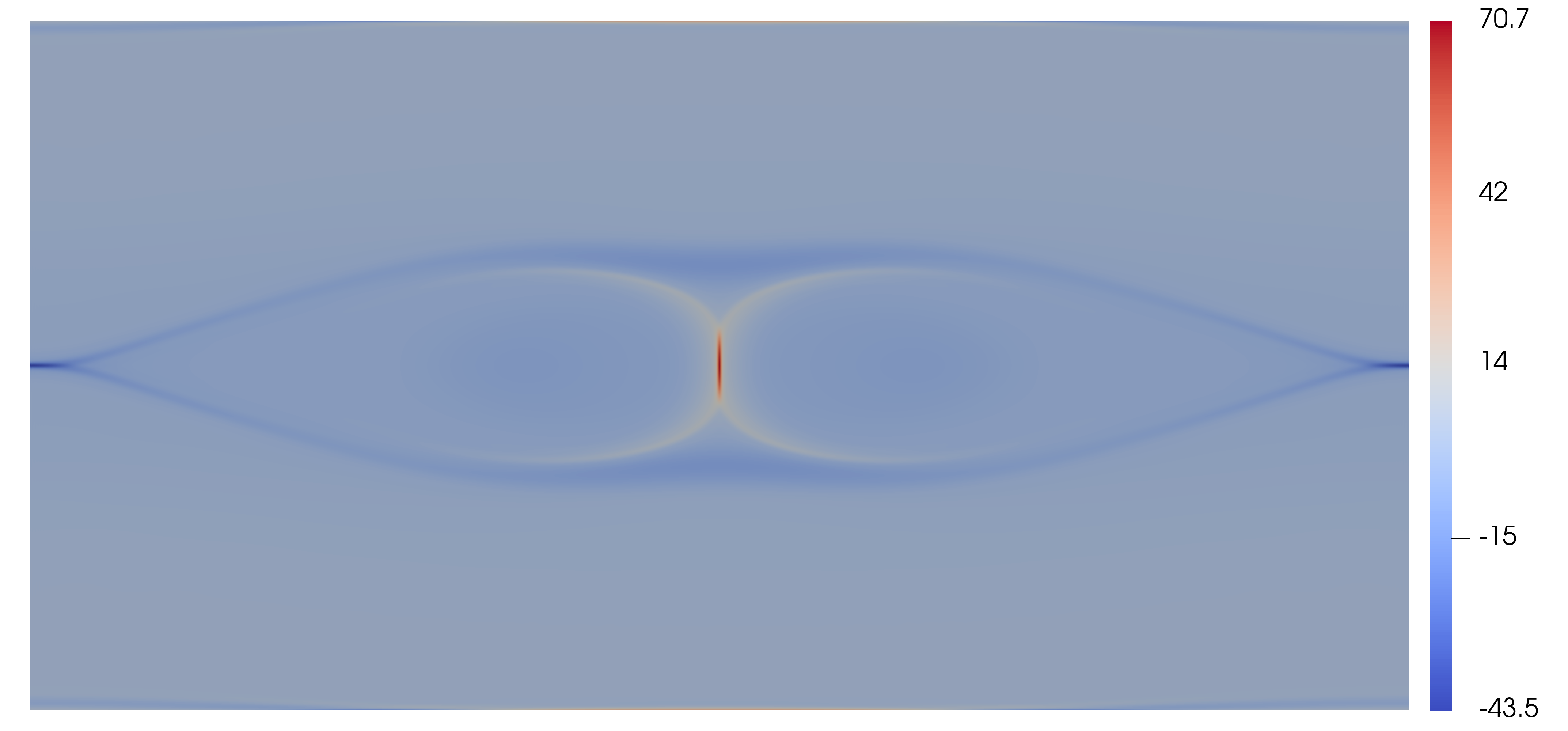}
        \caption{$t=2.5$}
    \end{subfigure}
    \caption{Island coalescence: current density $j_h$ at $t=0$, $1.5$, $2$ and $2.5$.}
    \label{fig:island}
\end{figure}

\begin{figure}[htbp]
    \centering
    \begin{subfigure}[b]{0.48\textwidth}
        \centering
        \includegraphics[width=\linewidth]{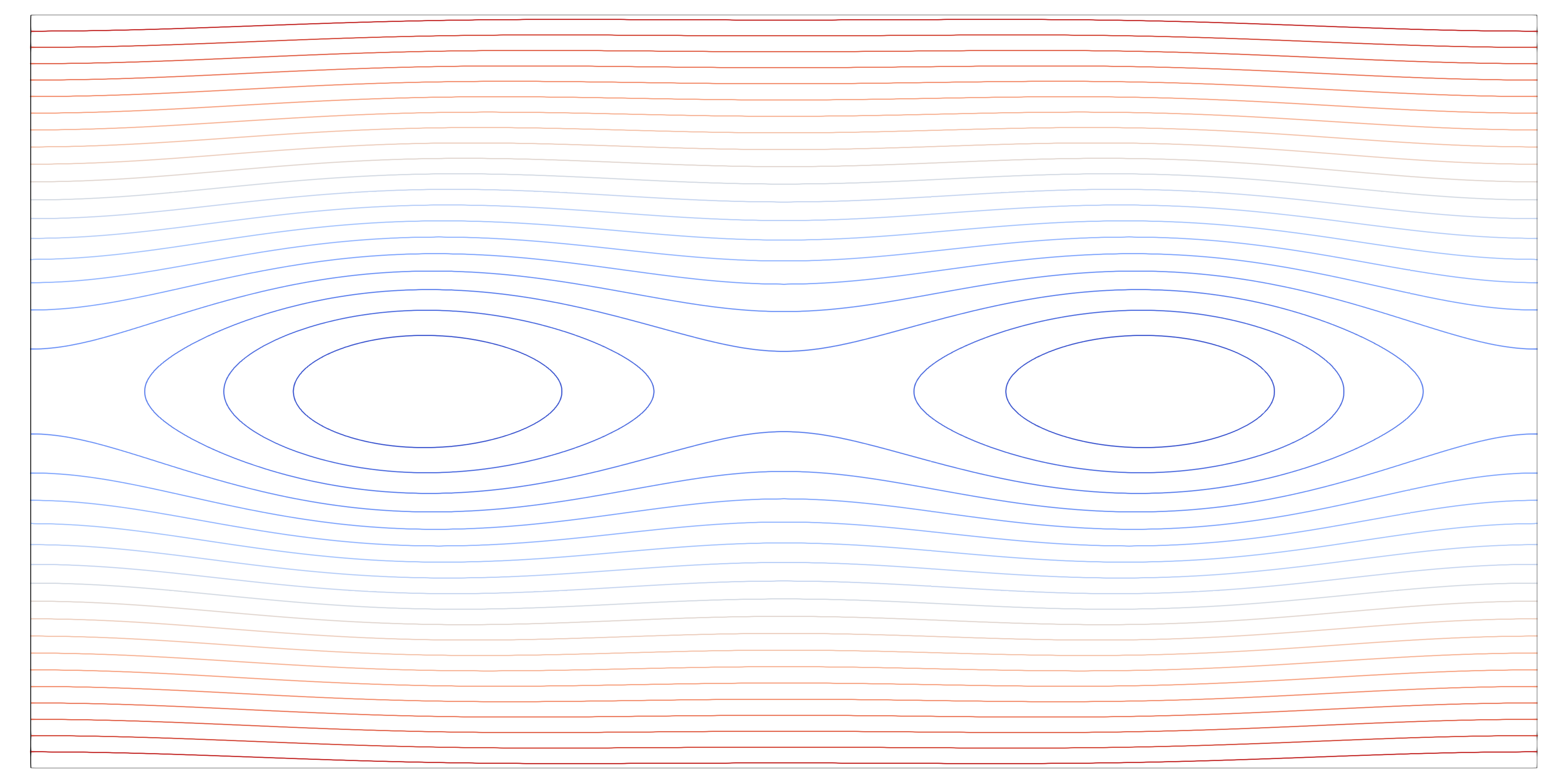}
        \caption{$t=0$}
    \end{subfigure}\hfill
    \begin{subfigure}[b]{0.48\textwidth}
        \centering
        \includegraphics[width=\linewidth]{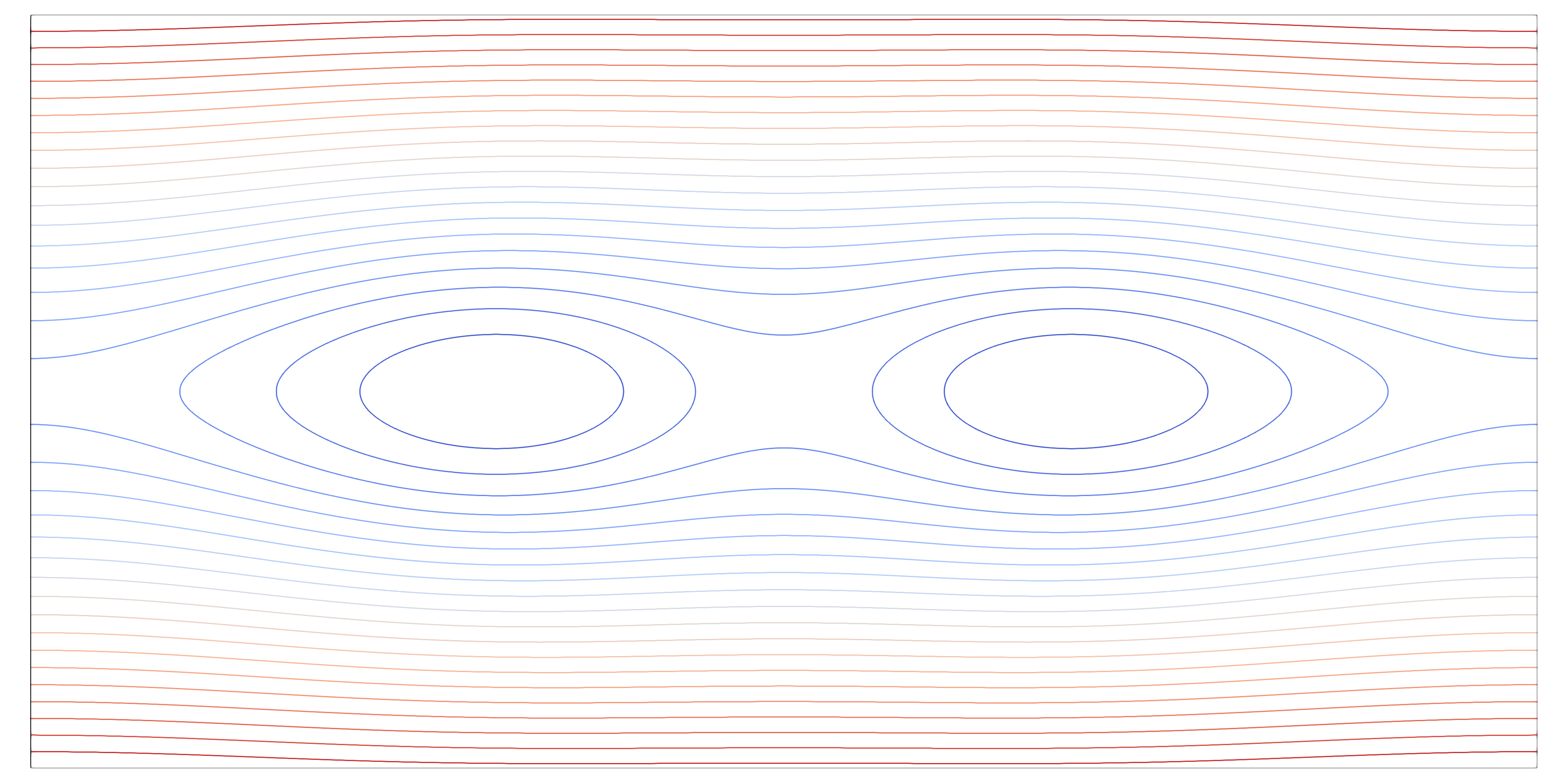}
        \caption{$t=1.5$}
    \end{subfigure}\\[1ex]
    \begin{subfigure}[b]{0.48\textwidth}
        \centering
        \includegraphics[width=\linewidth]{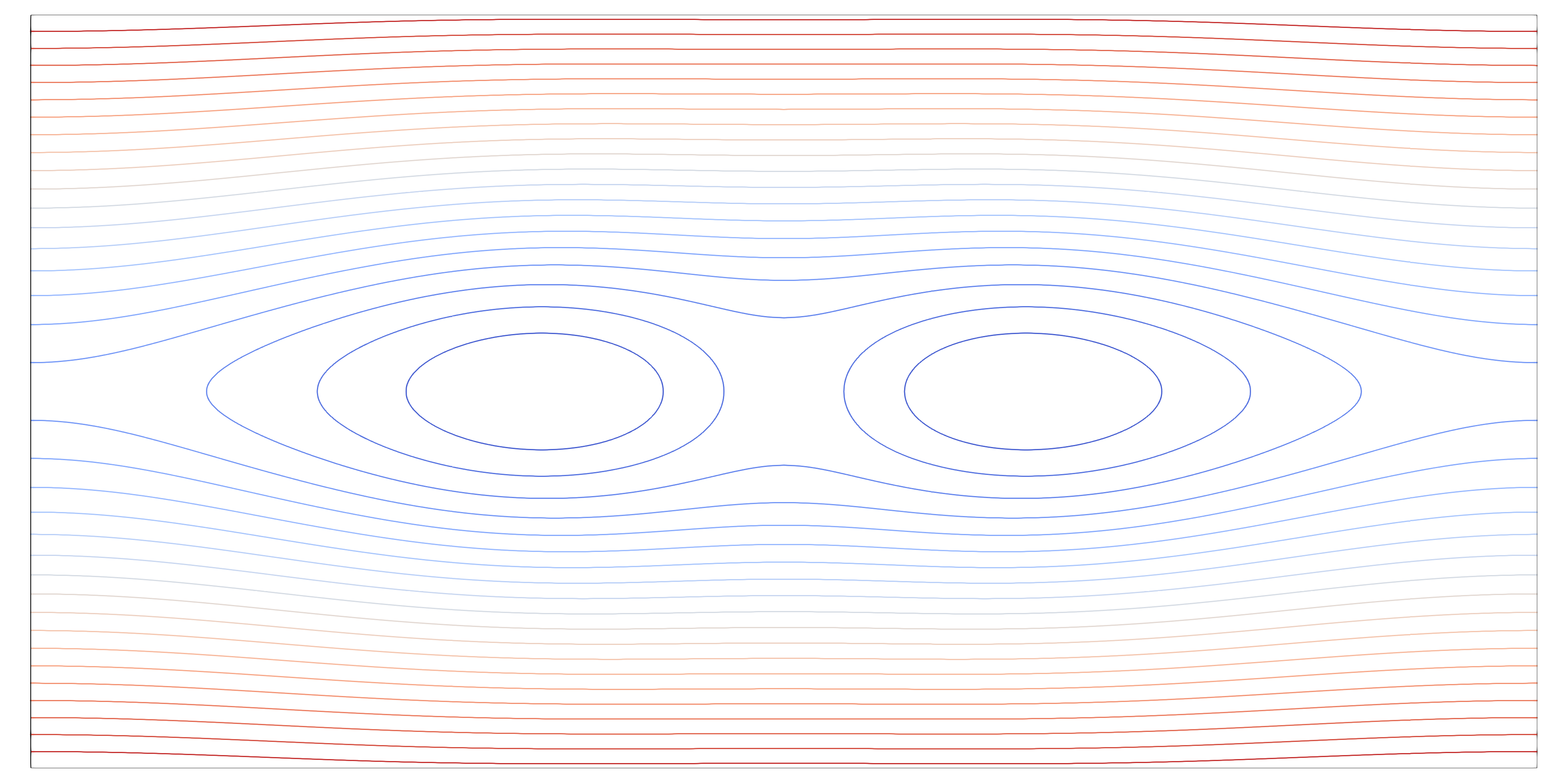}
        \caption{$t=2$}
    \end{subfigure}\hfill
    \begin{subfigure}[b]{0.48\textwidth}
        \centering
        \includegraphics[width=\linewidth]{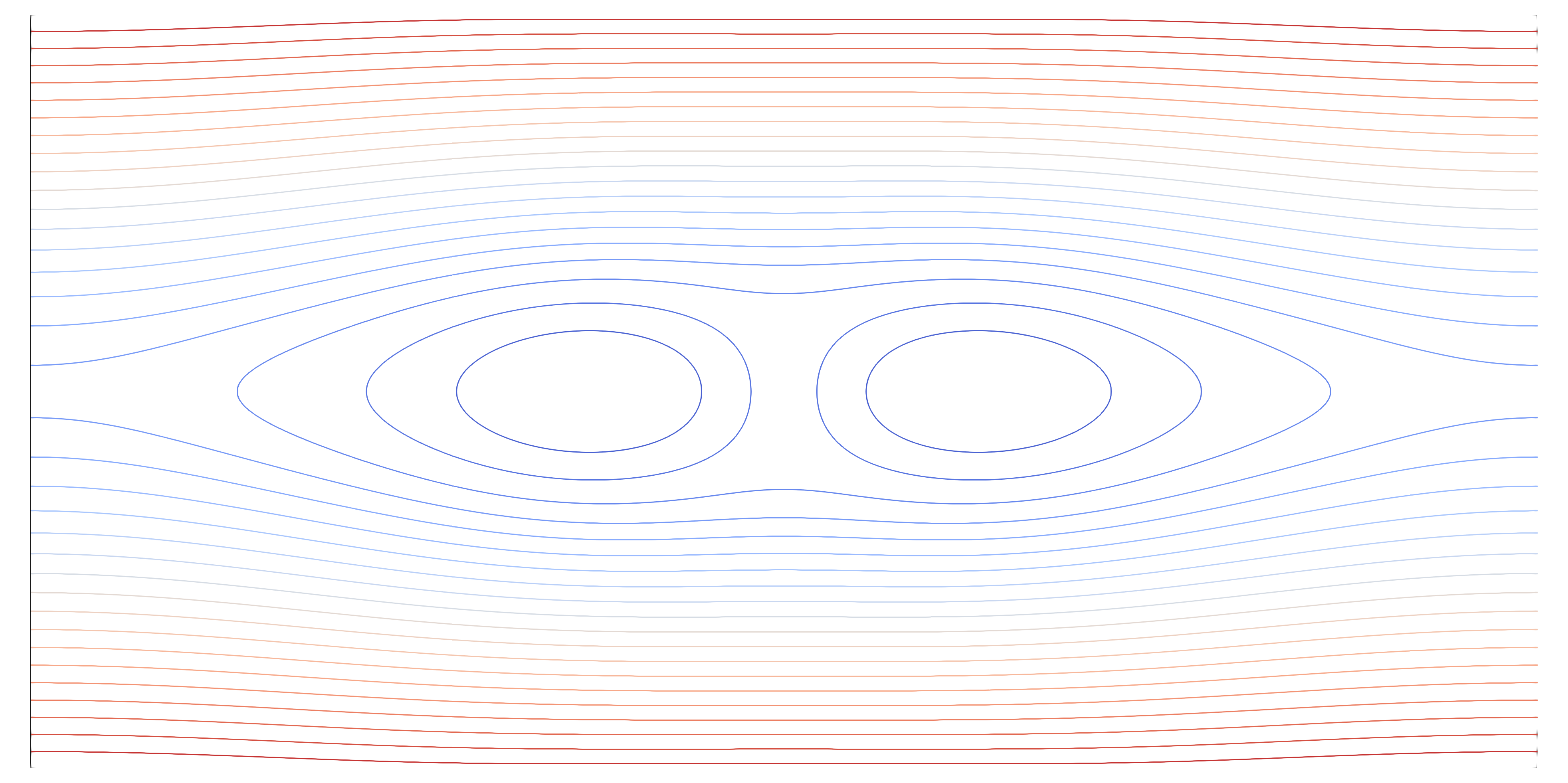}
        \caption{$t=2.5$}
    \end{subfigure}
    \caption{Island coalescence: magnetic potential $A_h$ (contours) at $t=0$, $1.5$, $2$ and $2.5$.}
    \label{fig:island-A}
\end{figure}

For the multigrid hierarchy we use a coarse mesh of $32\times 16$ cells and $5$ levels of refinement, resulting in a $1024\times 512$ mesh with $6{,}300{,}676$ degrees of freedom. We fix the time step $\Delta t = 0.005$ and integrate until $T=2.5$. The filter width is $\alpha=2h$. We use the solver tolerances of \Cref{sec:solver-development}; since the initial residual of each time step is small here, we additionally stop the nonlinear solver at an absolute residual of $10^{-8}$, the level reached by round-off in this computation. \Cref{fig:island,fig:island-A} report snapshots at $t=\{0, 1.5, 2, 2.5\}$.


At $t=0$ the two islands are separated by an X-point at $x=1$. The Lorentz force pulls them together, and by $t=1.5$ they have moved toward the center and a current sheet forms at the X-point. At $t=2.5$ the islands are pressed against each other and reconnect through a thin X-shaped current layer, while flux is also driven out along the side walls. These fine structures are resolved without spurious oscillations.

We further report the robustness of the solver with respect to the fluid and magnetic Reynolds numbers. We run the experiment until $T=2.5$ for six combinations of $\mathrm{Re},\mathrm{Re}_m\in\{10^2,10^4,10^6,10^8\}$ and report in \Cref{tab:reynolds-robustness} the average numbers of Newton steps and outer FGMRES iterations per Newton step, over all $500$ time steps and over $t\in[2,2.5]$, when the current sheet is strongest. Every time step takes two Newton steps, and the number of outer iterations stays between $1.46$ and $2$ per Newton step; it grows only slightly in the last phase for $\mathrm{Re}=\mathrm{Re}_m\ge10^6$, where the peak current exceeds $100$.


\begin{table}[htbp]
\centering
\caption{Robustness with respect to the Reynolds numbers for island coalescence up to $T=2.5$: average number of Newton steps (in parentheses) and average number of outer FGMRES iterations per Newton step, over all $500$ time steps and over $t\in[2,2.5]$, and the peak current $\|j_h\|_{L^\infty}$ at $T=2.5$.}
\label{tab:reynolds-robustness}

\begin{tabular}{rrccr}
\toprule
$\mathrm{Re}$ & $\mathrm{Re}_m$ & $t\in[0,2.5]$ & $t\in[2,2.5]$ & $\|j_h\|_{L^\infty}$ \\
\midrule
$10^{2}$ & $10^{2}$ & (2.00) 1.50 & (2.00) 1.50 & 3.4 \\
$10^{2}$ & $10^{8}$ & (2.00) 1.46 & (2.00) 1.50 & 80.4 \\
$10^{4}$ & $10^{4}$ & (2.00) 1.50 & (2.00) 1.50 & 32.4 \\
$10^{6}$ & $10^{6}$ & (2.00) 1.59 & (2.00) 1.92 & 115.1 \\
$10^{8}$ & $10^{2}$ & (2.00) 1.50 & (2.00) 1.50 & 3.5 \\
$10^{8}$ & $10^{8}$ & (2.00) 1.61 & (2.00) 2.00 & 129.1 \\
\bottomrule
\end{tabular}

\end{table}

\begin{figure}[htbp]
    \centering
    \includegraphics[width=\linewidth]{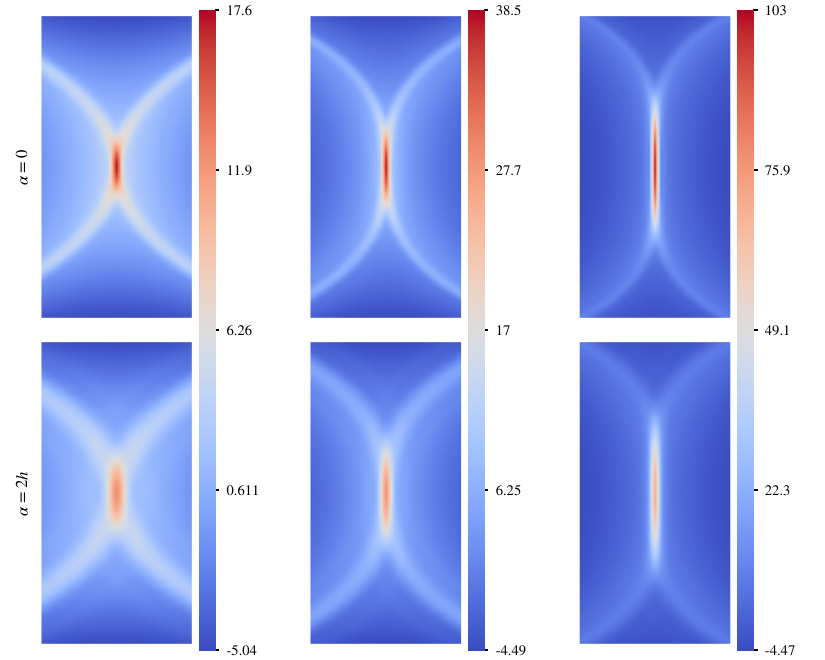}
    \caption{Island coalescence: current density $j_h$ on $[0.95,1.05]\times[-0.1,0.1]$ for $\alpha=0$ (top)
    and $\alpha=2h$ (bottom) at $t=2$, $2.25$ and $2.5$ (left to right). The two runs at each time share one
    color bar.}
    \label{fig:island-sheet-images}
\end{figure}

\begin{figure}[htbp]
    \centering
    \includegraphics[width=\linewidth]{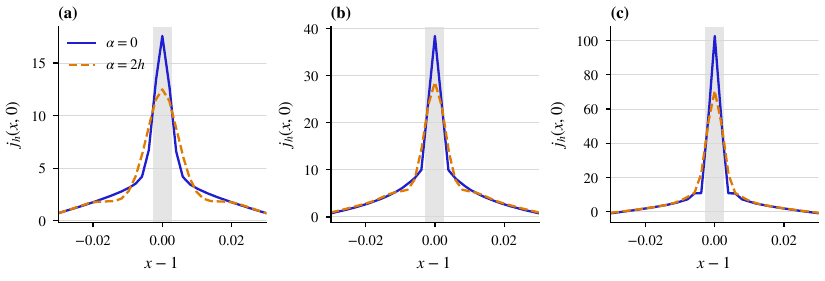}
    \caption{Island coalescence: current density across the sheet, $j_h(x,0)$, at (a) $t=2$, (b) $t=2.25$ and
    (c) $t=2.5$. The gray band has width $\alpha=2h$.}
    \label{fig:island-sheet-peak}
\end{figure}

\paragraph{Effect of the filter on the current sheet}
\Cref{fig:island-sheet-images,fig:island-sheet-peak} compare the current sheet for $\alpha=0$ and $\alpha=2h$ on
the same mesh. The filter smooths the sheet: with $\alpha=2h$ it is wider and its peak is lower,
while the field away from the sheet is unchanged. We thus observe that the regularization effect also applies to the small scales of the magnetic field, even though the magnetic field was not filtered directly.

\section{Conclusion}\label{sec:turbulence-conclusion}
We have constructed a structure-preserving finite element method for the two-dimensional
MHD-$\alpha$ model. Built on a discrete de~Rham complex, a discrete Helmholtz filter and the
implicit midpoint rule, the scheme preserves mass conservation and the magnetic Gauss law
pointwise and satisfies discrete dissipation laws for the energy, the cross helicity and the
mean square magnetic potential, which reduce to exact conservation in the ideal limit. The
discrete systems are solved by an augmented Lagrangian preconditioner with monolithic
geometric multigrid, whose iteration counts are observed to be robust with respect to the filter
width and the Reynolds numbers, and to the mesh size for small time steps. The scheme converges at the optimal rates, agrees with a pseudo-spectral
direct simulation at the large scales, and reproduces selective decay over long times.
Since the scheme introduces no numerical dissipation, the energy spectra measure the effect of
the filter alone. The filter leaves the large scales unchanged and suppresses small-scale velocity activity,
and when the cascade reaches the grid scale it reduces the energy accumulation there without
introducing additional dissipation.

A natural extension is the three-dimensional MHD-$\alpha$ model, where the mean square potential is replaced by the magnetic helicity, together with a helicity-preserving discretization \cite{he2025helicity,parkerE3,huHelicityconservativeFiniteElement2021,maoIncompressibilityDivB0Preserving2025,li2026linear}.


\section*{Acknowledgements}
This work was funded by the Engineering and Physical Sciences Research Council [grant number EP/W026163/1],
the Science and Technology Facilities Council [grant number UKRI/ST/B000495/1],
the Donatio Universitatis Carolinae Chair ``Mathematical modelling of multicomponent systems'',
the UKRI Digital Research Infrastructure Programme through the Science and Technology Facilities Council's Computational Science Centre for Research Communities (CoSeC),
the European Research Council (ERC Starting Grant, project 101164551 GeoFEM), and
a Royal Society University Research Fellowship (URF$\backslash$R1$\backslash$221398).

\section*{Code availability}
The scheme is implemented using Firedrake \cite{FiredrakeUserManual} and PETSc \cite{petsc-user-ref}. The code and the scripts that generate all tables and figures are archived on Zenodo \cite{2dmhd-alpha}. The reference DNS in \Cref{fig:2d-ozt-spectrum-nu1e-4} was computed with Dedalus \cite{burns2020dedalus}.

\section*{Declaration of generative AI and AI-assisted technologies in the writing process}
During the preparation of this work the authors used Claude Sonnet 5.5, Opus 5.5, Fable 5.2 (Anthropic) and GPT Sol 5.6 (OpenAI) in order to improve the clarity, readability, and language of the manuscript and to assist with \LaTeX{} typesetting. After using these tools, the authors reviewed and edited the content as needed and take full responsibility for the content of the publication. The AI tools were also used to check derivations and to help debug code; all mathematical statements, proofs and numerical results were produced and verified by the authors.

\appendix
\section{Computation of the energy spectra}\label{app:spectra}

On $\Omega=(0,2\pi)^2$, let
$\hat{\boldsymbol v}_{\boldsymbol\kappa}=(2\pi)^{-2}\int_\Omega\boldsymbol v\,
e^{-\mathrm i\boldsymbol\kappa\cdot\boldsymbol x}\,\mathrm dx$. The kinetic energy
spectrum is the shell sum
\begin{equation}\label{eq:shell}
  E^\alpha_K(k)=\sum_{k\le|\boldsymbol\kappa|<k+1}\tfrac12\operatorname{Re}
  \big(\hat{\boldsymbol u}_{\boldsymbol\kappa}\cdot
  \overline{\hat{\overline{\boldsymbol u}}_{\boldsymbol\kappa}}\big),
\end{equation}
so that $(2\pi)^2\sum_{k\ge0}E^\alpha_K(k)=\frac12(\boldsymbol u,\overline{\boldsymbol u})$.
We compute this spectrum by the FFT of point samples of the finite element fields. The sampling
points are chosen so that the computed spectrum sums exactly to the assembled energy
$\frac12(\boldsymbol u_h,\overline{\boldsymbol u}_h)$, and the sampling grid is fine enough that
aliasing does not affect the plotted range $k\le N/2$.

\begin{figure}[htbp]
\centering
\begin{tikzpicture}[
  font=\small,
  lo/.style={circle,fill=black,inner sep=0pt,minimum size=3.2pt},
  up/.style={circle,draw=black,fill=white,inner sep=0pt,minimum size=3.2pt,line width=0.5pt},
  vec/.style={-latex,line width=0.7pt}]
\begin{scope}[scale=2.4]
  \foreach \i in {0,1} \foreach \j in {0,1}
    \fill[gray!12] (\i,\j) -- (\i+1,\j) -- (\i,\j+1) -- cycle;
  \draw[line width=0.6pt] (0,0) grid (2,2);
  \foreach \i in {0,1} \foreach \j in {0,1}
    \draw[line width=0.6pt] (\i,\j+1) -- (\i+1,\j);
  \foreach \i in {0,...,7} \foreach \j in {0,...,7} {
    \pgfmathsetmacro{\x}{(\i+0.31)/4}
    \pgfmathsetmacro{\y}{(\j+0.47)/4}
    \pgfmathsetmacro{\lx}{\x-floor(\x)}
    \pgfmathsetmacro{\ly}{\y-floor(\y)}
    \pgfmathparse{\lx+\ly>1 ? 1 : 0}
    \ifnum\pgfmathresult=1 \node[up] at (\x,\y) {}; \else \node[lo] at (\x,\y) {}; \fi
  }
  \node at (1,2.25) {(a)\ $M=4N$, $\boldsymbol o=(0.31,0.47)$};
\end{scope}
\begin{scope}[xshift=6.2cm,scale=4.2]
  \fill[gray!12] (0,0) -- (1,0) -- (0,1) -- cycle;
  \draw[line width=0.6pt] (0,0) rectangle (1,1);
  \draw[line width=1.1pt] (0,1) -- (1,0);
  \foreach \p in {(0.12,0.12),(0.42,0.10),(0.10,0.45)}
    \draw[vec] \p -- ++(0.16,0.032);
  \foreach \p in {(0.86,0.58),(0.55,0.72),(0.80,0.80)}
    \draw[vec] \p -- ++(0.032,0.16);
  \node at (0.30,0.30) {$\boldsymbol u_{T_1}$};
  \node at (0.70,0.60) {$\boldsymbol u_{T_2}$};
  \draw[red!75!black,line width=1.1pt] (0.47,0.47) -- (0.53,0.53) (0.47,0.53) -- (0.53,0.47);
  \node[red!75!black,anchor=north,align=center,font=\footnotesize] at (0.5,-0.05)
    {$M=N$, $\boldsymbol o=(\tfrac12,\tfrac12)$: $\boldsymbol u_{T_1}\neq\boldsymbol u_{T_2}$};
  \node at (0.5,1.18) {(b)\ one cell, $\boldsymbol u_h\in\bm Z_h$};
\end{scope}
\end{tikzpicture}
\caption{Sampling of the Raviart--Thomas fields. (a) Sample points on $2\times2$ cells of
the mesh; no sample lies on a side or diagonal of a cell. (b) The field is constant on each
triangle, so a sample on the diagonal takes two different values.}
\label{fig:sampling}
\end{figure}
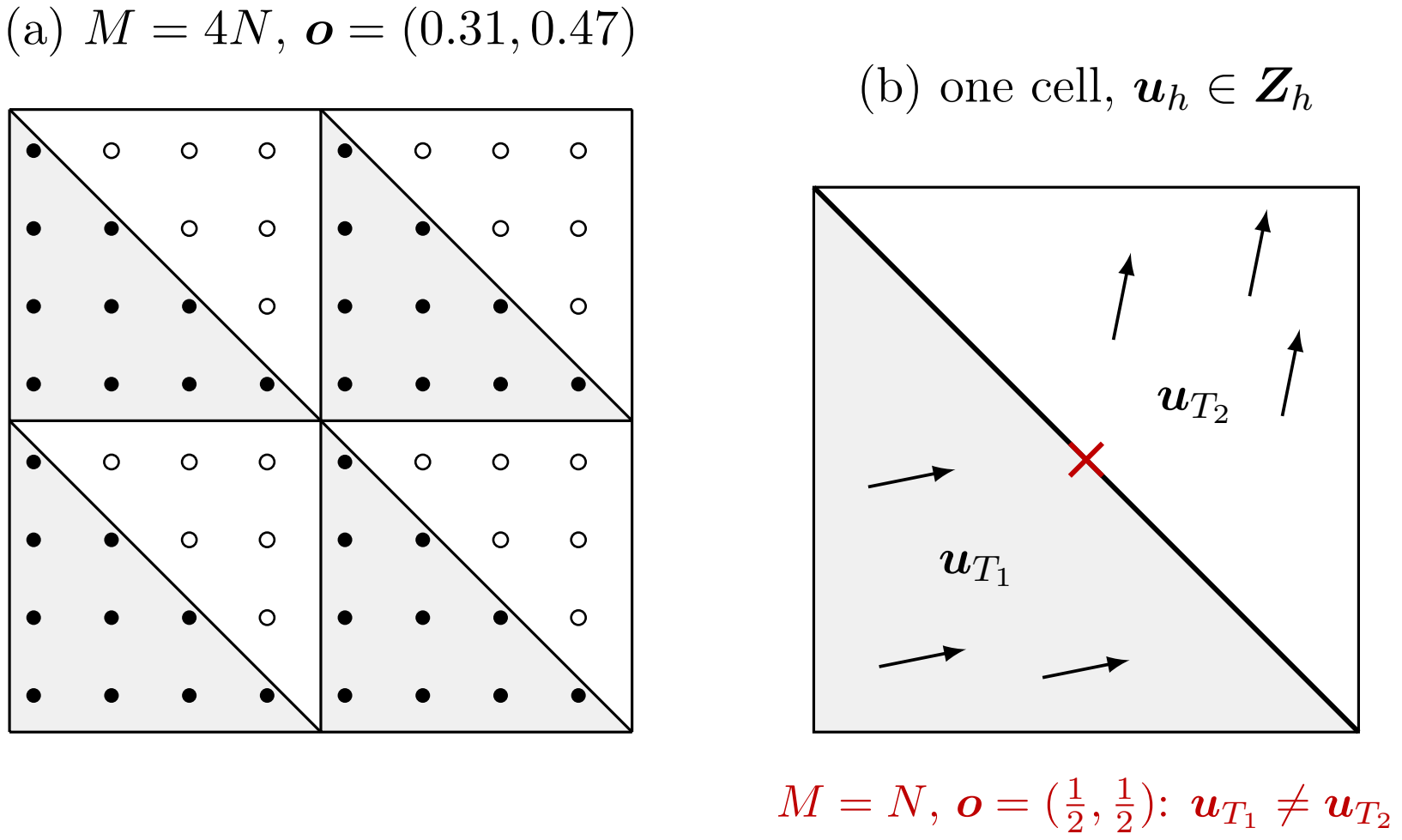

\textit{Sampling points.} When sampling a finite element field to compute the energy spectrum, the samples must avoid being taken on the edges. For $k=1$, the discrete fields are piecewise constant: a lowest-order Raviart--Thomas
function is $\boldsymbol a_T+b_T\boldsymbol x$ on each cell, with divergence $2b_T$,
and $\boldsymbol u_h$, $\overline{\boldsymbol u}_h$ are divergence-free pointwise by
\Cref{thm:2d-sp}(a). A point value anywhere inside a cell therefore equals the cell value.
On an edge, however, the tangential component jumps and the field takes two values
(\Cref{fig:sampling}(b)). We therefore sample at
$\boldsymbol x_{ij}=\frac{2\pi}{M}(i+o_1,j+o_2)$, $i,j=0,\dots,M-1$, with $M=rN$,
$r\in\mathbb N$; the offsets $\boldsymbol o=(0.31,0.47)$ satisfy
$o_1,o_2,o_1\pm o_2\notin\mathbb Z$, which keeps every sample off the sides and
diagonals of the cells (\Cref{fig:sampling}(a)).

\textit{Sampled versus assembled energy.} We show that the sampled energy then equals the assembled energy
$\frac12(\boldsymbol u_h,\overline{\boldsymbol u}_h)$. The samples are not distributed in
proportion to the area, since for $r=4$ the two triangles of a square receive $10$ and $6$
samples. The equality holds because every square is split by the same diagonal, from its
upper-left to its lower-right corner. For $\psi_h\in V_h$ the
$x$-derivative on the upper triangle of a square equals that on the lower triangle of the
square above it, and the $y$-derivative equals that on the lower triangle of the square to
its right. By \eqref{eq:exact-2d} we may write $\boldsymbol u_h=\curl\psi_h+\boldsymbol c$ and
$\overline{\boldsymbol u}_h=\curl\overline\psi_h+\overline{\boldsymbol c}$ with
$\psi_h,\overline\psi_h\in V_h$ and constant vectors $\boldsymbol c,\overline{\boldsymbol c}$.
Then $\boldsymbol u_T\cdot\overline{\boldsymbol u}_T$ is a sum of products of two derivatives
of $\psi_h$ and $\overline\psi_h$ in the same direction, of single such derivatives times
constants, and of the constant $\boldsymbol c\cdot\overline{\boldsymbol c}$. Each of these
terms is shifted from the upper triangles to the lower triangles as above, so the sums of
$\boldsymbol u_T\cdot\overline{\boldsymbol u}_T$ over the lower and over the upper triangles
coincide, and since the samples are split in
the same way in every square,
$(\frac{2\pi}{M})^2\sum_{i,j}\boldsymbol u_h(\boldsymbol x_{ij})\cdot
\overline{\boldsymbol u}_h(\boldsymbol x_{ij})=(\boldsymbol u_h,\overline{\boldsymbol u}_h)$.

\textit{Discrete spectrum.} In \eqref{eq:shell} we replace the Fourier coefficients by the discrete coefficients
$\hat{\boldsymbol v}^M_{\boldsymbol\kappa}=M^{-2}\sum_{i,j}\boldsymbol v(\boldsymbol x_{ij})
e^{-\mathrm i\boldsymbol\kappa\cdot\boldsymbol x_{ij}}$ of the samples, computed with the FFT of
NumPy \cite{harris2020array}. Here $M$ is even, $\boldsymbol\kappa$ is restricted to the representative index set $\mathcal K_M:=\{-M/2,\ldots,M/2-1\}^2$, and the shell sums in \eqref{eq:shell} are restricted to $\mathcal K_M$. Parseval's identity holds exactly for these coefficients, so the
spectrum sums to the sampled energy, which equals the assembled energy. The constants $\boldsymbol c,\overline{\boldsymbol c}$ enter only the mode $\boldsymbol\kappa=0$.
\Cref{tab:parseval} shows the check for an Orszag--Tang computation with $\nu=\eta=10^{-3}$ on
the $256\times256$ mesh, where the relative differences are at the level of round-off. For the
computations of \Cref{fig:2d-ozt-spectrum-nu1e-4} we observed relative differences below
$3\times10^{-12}$ for $M=N$, $2N$ and $4N$, which is the level of the remaining discrete divergence
when the nonlinear systems are solved iteratively. We use $M=4N$ and show the spectra for
$k\le N/2$. For the reference we evaluate \eqref{eq:shell} with $\overline{\boldsymbol u}=\boldsymbol u$ and the
Fourier coefficients of the pseudo-spectral solution, obtained by the FFT of its vorticity on the
$1024\times1024$ collocation grid, $|\hat{\boldsymbol u}_{\boldsymbol\kappa}|^2
=|\hat w_{\boldsymbol\kappa}|^2/|\boldsymbol\kappa|^2$; since this solution is a trigonometric
polynomial, these coefficients are exact and no sampling is needed.

\textit{Aliasing.} The discrete coefficients differ from the exact ones by aliasing,
$\hat{\boldsymbol v}^M_{\boldsymbol\kappa}=\sum_{\boldsymbol m\in\mathbb Z^2}
e^{2\pi\mathrm i\,\boldsymbol m\cdot\boldsymbol o}\,\hat{\boldsymbol v}_{\boldsymbol\kappa+M\boldsymbol m}$.
For $k\le N/2$ only wavenumbers $|\boldsymbol\kappa+M\boldsymbol m|\ge 3.5N$, corresponding to length scales far below the mesh
scale, are folded back. Increasing $M$ beyond $4N$ leaves the spectra of
\Cref{fig:2d-ozt-spectrum-nu1e-4} unchanged.
Aliasing therefore does not affect the plotted spectra.

\textit{Mean square potential.} The spectrum of the mean square potential in \Cref{fig:A2-inverse-cascade} is computed in the same
way from point samples of $A_h\in V_h$, $A(k)=\sum_{k\le|\boldsymbol\kappa|<k+1}|\hat A^M_{\boldsymbol\kappa}|^2$,
with $M=4N$; since $A_h$ is continuous, the samples need not avoid the edges. At each time the
spectrum is normalized by its maximum $A_{\max}=\max_k A(k)$, so that only its shape is compared.

\begin{table}[htbp]
\centering
\caption{Sampled versus assembled kinetic energy for the Orszag--Tang test with $\nu=\eta=10^{-3}$
on the $256\times256$ mesh: assembled energy $\frac12(\bm u_h,\ov{\bm u}_h)$ and relative difference
$|(2\pi)^2\sum_kE^\alpha_K(k)-\frac12(\bm u_h,\ov{\bm u}_h)|/\frac12(\bm u_h,\ov{\bm u}_h)$ for the
sampling grids $M=N$ and $M=4N$.}
\label{tab:parseval}
\small
\begin{tabular}{ccccc}
\toprule
$\alpha$ & $t$ & assembled & $M=N$ & $M=4N$ \\
\midrule
$0$  & $0.5$ & 21.694511 & $3.1\times10^{-15}$ & $3.3\times10^{-16}$ \\
     & $1.5$ & 17.430970 & $4.1\times10^{-16}$ & $6.1\times10^{-16}$ \\
     & $2$   & 15.372711 & $4.6\times10^{-16}$ & $2.3\times10^{-16}$ \\
\midrule
$2h$ & $0.5$ & 21.505903 & $1.7\times10^{-16}$ & $9.9\times10^{-16}$ \\
     & $1.5$ & 16.455748 & $2.2\times10^{-16}$ & $8.6\times10^{-16}$ \\
     & $2$   & 14.376904 & $2.2\times10^{-15}$ & $2.0\times10^{-15}$ \\
\midrule
$4h$ & $0.5$ & 20.989329 & $1.0\times10^{-15}$ & $8.5\times10^{-16}$ \\
     & $1.5$ & 15.744024 & $2.3\times10^{-15}$ & $3.9\times10^{-15}$ \\
     & $2$   & 13.655930 & $2.1\times10^{-15}$ & $1.3\times10^{-15}$ \\
\bottomrule
\end{tabular}
\end{table}

\bibliographystyle{elsarticle-harv} 
\bibliography{ref}
\end{document}